\documentclass[12pt]{amsart}
\usepackage{amssymb}
\usepackage{eucal}
\newcommand{\rr}{\mathbb{R}}
\newcommand{\zz}{\mathbb{Z}}
\newcommand{\cc}{\mathbb{C}}

\newtheorem{lemma}{Lemma}[section]
\newtheorem{proposition}{Proposition}[section]
\newtheorem{theorem}{Theorem}[section]
\newtheorem{definition}[theorem]{Definition}

\author{C. E. Kenig, G. Ponce, C. Rolvung, and L. Vega}
\title[ Schr\"odinger Equation]
{The General Quasilinear Ultrahyperbolic Schr\"odinger Equation}

\begin{document}
\maketitle

\numberwithin{equation}{section}

\section{Introduction}

In this article we consider nonlinear Schr\"odinger equations of the form

\begin{equation}
\begin{cases}
\label{1.1}
\begin{aligned}
&\partial_t u=-\, i\,\partial_{x_j}(a_{jk}(x,t,u,\bar u,\nabla u,\nabla \bar u) \partial_{x_k}u)\\
&+\vec b_1(x,t,u,\bar u, \nabla u,\nabla \bar u)\cdot \nabla u 
+\vec b_2(x,t,u,\bar u, \nabla u,\nabla \bar u)\cdot \nabla \bar u\\
&+ c_1(x,t, u,\bar u)u + c_2(x,t, u,\bar u) \bar u +f(x,t),
\end{aligned}
\end{cases}
\end{equation}
where $x\in\rr^n$, $t>0$, and $A=(a_{jk}(\cdot))_{j,k=1,..,n}$ is a real, symmetric matrix.

Our aim is to study the existence, uniqueness and regularity of local solutions to
the initial value problem (IVP) associated to the equation (1.1).

In the case where $A=(a_{jk}(\cdot))_{j,k=1,..,n}$ is assumed to be elliptic the local solvability of the
IVP associated to (1.1) was recently established in \cite{KePoVe4}. Hence, in this work we should be concerned with the case
where $(a_{jk}(\cdot))_{j,k=1,..,n}$ is just a non-degenerate matrix.

Equations of the form described in (1.1) with  $A=(a_{jk}(\cdot))_{j,k=1,..,n}$ merely invertible 
arise  in  water wave problems,
and as higher dimensions completely integrable models,
see for example \cite {AbHa}, \cite{DaSt}, \cite{DjRe}, \cite{Is}, \cite{s-s}, and \cite{ZaSc}.

There are significant  differences in the arguments required for the local solvability in the case where 
$A$  is a non-degenerate matrix in comparison with
the elliptic case treated in \cite{KePoVe4}. To illustrate them as well as 
to review some of the previous related results we consider first 
the semi-linear equation with constant coefficients
(for more details and further references and comments see \cite{KePoVe3}, 
\cite{KePoVe4}, 
\cite{KePoRoVe}, and references therein)
\begin{equation}
\label{1.2a}
\partial_tu=-i(\partial_{x_1}^2+..+\partial_{x_k}^2-\partial_{x_{k+1}}^2-..-\partial_{x_n}^2)u+P(u,\bar 
u,\nabla_xu, \nabla_x\bar 
u),
\end{equation}
where $P(\cdot)$ is the non-linearity (by simplicity a polynomial in its variables without constant or linear terms).

In \cite{KePoVe2} based on the smoothing effects (homogeneous and inhomogeneous,
see \cite{Kt}, \cite{KrFa}, \cite{c-s}, \cite{s}, \cite{v}, \cite{KePoVe1}) associated to the group
$\{e^{-it(\partial_{x_1}^2+..+\partial_{x_k}^2-\partial_{x_{k+1}}^2-..-\partial_{x_n}^2)}\,:\,t\in\rr\}$
the local wellposedness for \lq\lq small" data for the IVP associated to \eqref{1.2a} was deduced.
In \cite{HaOz} for the one dimensional case ($n=1$), Hayashi and Ozawa eliminated the size restriction on the 
data
in \cite{KePoVe2}. Their argument was based on a change of variable which transforms the equation into
a new system without the term $\partial_xu$, so that the standard energy estimate yields the
desired local result. In \cite{Ch} Chihara, for the elliptic case (i.e. $k=n$ in \eqref{1.2a}),
removed the  size restriction on the data in any dimension $n$. Roughly speaking, the argument there first 
uses the ellipticity
to diagonalize the  system for $(u,\bar u)$, and then introduced an operator $K$ so that the commutator
$i[K;\Delta]$ \lq\lq controls" the term $K\vec b(x)\cdot\nabla_x$. This is achieved by combining  some result of 
Doi  \cite{Do1}
concerning the local smoothing effects in the solution with  the sharp version of Garding inequality.

 If instead of \lq\lq controlling" it one asks for the operator $K$ to verify that
\begin{equation}
\label{commutator}
-i[K;\Delta]+K\vec b(x)\cdot \nabla_x = 0 \,+\;\text{order zero},
\end{equation}
one finds that $K$ has symbol
\begin{equation}
\label{kks}
k(x,\xi)=exp(-\int_0^{\infty}\vec b(x+2s\xi)\cdot \xi ds),
\end{equation}
which is in the non-standard class studied by Craig-Kappeler-Strauss \cite{CrKaSt}. In particular, it satisfies that 
$$
|\partial_x^{\alpha}\partial_{\xi}^{\beta} k(x,\xi)|
\leq c_{\alpha \beta}\langle x\rangle^{|\beta|}\langle \xi\rangle^{-|\beta|},\;\;\;\;\;\;\alpha,\,\beta\in(\zz^+)^n.
$$ 
 
However, in the non-elliptic case with coefficients depending just on the space variable $x$,
the geometric assumptions in \cite{CrKaSt} (Chapter 3.1, section 3.1) does not
hold for the relevant symbols.
 
 The local wellposedness of the IVP associated to the equation in \eqref{1.2a} ($1<k\leq n$)
was established in \cite{KePoVe3}. 
The method of proof there, among other arguments, utilizes  the 
 symbol class $\;S^0_{0,0}\;$ of Calder\'on-Vaillancourt \cite{CaVa}.
However,  this approach does not seem  to extend
to the variable coefficients case.

 Next, we consider the linear IVP  
\begin{equation}
\label{1.2linear}
\begin{cases}
\begin{aligned}
&\partial_tu= - i  \partial_{x_k}a_{jk}(x)\partial_{x_j} u+\vec b_1(x)\cdot \nabla u+
\vec b_2(x)\cdot \nabla\overline u +f(x,t),\\
& u(x,0)=u_0(x).
\end{aligned}
\end{cases}
\end{equation}

 We recall the notion of the bicharacteristic flow associated to the symbol
of the principal part of the operator 
$\,-\partial_{x_k}a_{jk}(x)\partial_{x_j}$, i.e. $h(x,\xi)=a_{jk}(x)\xi_k\xi_j$.

Let $(X(s;x_0,\xi_0),\Xi(s;x_0,\xi_0))$ denote the solution of the Hamiltonian system
\begin{equation}
\label{bic}
\begin{cases}
\begin{aligned}
&\frac{d}{ds} X_j(s;x_0,\xi_0)=
2 a_{jk}(X(s;x_0,\xi_0))\,\Xi_k(s;x_0,\xi_0)=\partial_{\xi_j}h,\\
\\
&\frac{d}{ds} \Xi_j(s;x_0,\xi_0)=-\partial_{x_j}a_{lk}(X(s;x_0,\xi_0))
\Xi_k \Xi_l(s;x_0,\xi_0)=-\partial_{x_j}h,
\end{aligned}
\end{cases}
\end{equation}

for $j=1,..,n$, with data
$$
(X(0;x_0,\xi_0),\Xi(0;x_0,\xi_0))=(x_0,\xi_0).
$$

 Under mild regularity assumptions on the coefficients $a_{jk}(x)$'s the bicharacteristic flow 
$(X(s;x_0,\xi_0),\Xi(s;x_0,\xi_0))$
is defined in the interval $s\in(-\delta,\delta)$ with
$\delta=\delta(x_0,\xi_0)>0$ depending continuously on $(x_0,\xi_0)$.

 If the operator $ - a_{jk}(x)\partial^2_{x_kx_j}$  is elliptic, i.e.
$(a_{jk}(x))$ is positive definite, using that the flow preserves $h$, i.e.
$$
H_{h}h=\sum_{j=1}^n\,(\partial_{\xi_j}h \partial_{x_j}p-\partial_{x_j}h \partial_{\xi_j}p)=\{h,p\}=0, 
$$
one has that 
$$
\nu^{-2}|\xi_0|^2\leq |\Xi(s;x_0,\xi_0)|^2\leq \nu^2|\xi_0|^2.
$$
Hence $\delta=\infty$, i.e. the bicharacteristic flow is globally defined.
In the non-elliptic case one needs to prove it.
 
In \cite{Ic} Ichinose established that for the $L^2$-local
wellposedness of IVP \eqref{1.2linear}, with $-\partial_{x_k} a_{jk}(x)\partial_{x_j}$  elliptic
and $b_2(x)\equiv 0$, it is necessary that
\begin{equation}
\label{ichi}
\sup_{x\in\rr^n, \omega \in  S^{n-1}, R>0}\;
\left|\,Im \;\int_0^{R}
\vec b_1(X(s;x,\omega))\cdot\Xi(s;x,\omega) ds\right| < \infty.
\end{equation}

This extends previous results of Mizohata \cite{Mi} and Takeuchi \cite{Ta} for the constant coefficient case,
where $(X(s;x_0,\xi_0), \Xi(s;x_0,\xi_0))=(x_0+2s\xi_0,\xi_0)$.  
 Notice that in this variable coefficient case the \lq\lq integrating factor" in \eqref{kks} reads
\begin{equation}
\label{intfac}
k(x,\xi)=exp\,(-\int_0^{\infty} \vec b_1(X(s;x,\xi))\cdot\Xi(s;x,\xi) ds).
\end{equation}

 The condition \eqref{ichi} justifies the following \lq\lq non-trapping" assumption. The bicharacteristic flow 
\eqref{bic} is non-trapped if the set
$$
\{X(s;x_0,\xi_0)\,:\,s\in\rr^+\}
$$
is  unbounded in $\rr^n$ for each $(x_0,\xi_0)\in\rr^n\times\rr^n-\{0\}$.

 As it was already mentioned the IVP for the equation \eqref{1.1}
with $(a_{jk}(\cdot))_{j,k=1,..,n}$ elliptic was studied in \cite{KePoVe4}.
There the local solvability was obtained under 
regularity and decay  assumptions on the coefficients.
Also a non-trapping  character of the
bicharacteristic flow 
associated to the principal symbol of the elliptic operator 
$$
- \,\partial_{x_j}(a_{jk}(x,t,u,\bar u, \nabla_x u,\nabla_x \bar u)\partial_{x_k}
$$
when evaluated at the data $u(x,0)=u_0(x)$ and at a time $t=0$ was assumed, i.e. the bicharacteristic flow
associated to the  symbol 
\begin{equation}
\label{symbol}
h(u_0)=h_{u_0}(x,\xi)= a_{jk}(x,0, u_0, \bar u_0,\nabla_x u_0, \nabla_x\bar 
u_0)\xi_j\xi_k.
\end{equation}

The proofs of the  semi-linear results in \cite{KePoVe2}, \cite{KePoVe3}
follow a fixed point theorem (via contraction principle). This 
approach does not extend to the 
quasi-linear case, thus in \cite{KePoVe4} the proof 
used  the so called \lq\lq artificial viscosity" method. 

 We recall that one of the advantages of the contraction principle approach is
that it also provides the continuity (in fact, the analyticity) of the solution upon the data.
  
Returning to the non-degenerate case, in  \cite{KePoRoVe}  the following semi-linear IVP was studied
\begin{equation}
\label{1.2semi-l}
\begin{cases}
\begin{aligned} 
&\partial_tu= - i  \partial_{x_k}a_{jk}(x)\partial_{x_j} u+\vec b_1(x)\cdot \nabla u+
\vec b_2(x)\cdot \nabla\overline u\\
&\;\;\;\;+c_1(x) u+c_2(x)\overline u+
P(u,\nabla u,\overline u,\nabla \overline u),\\
& u(x,0)=u_0(x),
\end{aligned}
\end{cases}
\end{equation}
where the non-linearity $P$ is given by a  polynomial without linear or constant terms.
Under  the following assumptions:

(a) (Non-degeneracy) There exists $\nu\in (0,1)$
$$
\nu|\xi|\le |A(x)\xi|\le \nu^{-1} |\xi|,\;\;\;\;\,\,\,\,x,\xi\in\rr^n,
$$
where $A(x)=(a_{jk}(x))_{j,k=1,..,n}$.

(b)  (Asymptotic Flatness) There exist $c_0>0$ and $N,\, \tilde M\in\zz^+$ large enough such 
that 
\begin{equation}
\label{ddeca}
|\partial_x^{\alpha}(a_{jk}(x)-a^0_{jk})|\leq \frac{c_0}{\langle x\rangle^N},\;\;\;|\alpha|\leq 
\tilde M,\;\;\,j,k=1,..,n,
\end{equation}
where 
\begin{equation}
\label{matrix}
A_h=(a^0_{jk})_{j,k=1,..,n}=
\begin{pmatrix}I_{k\times k}&0\\0&-I_{(n-k)\times(n-k)}
\end{pmatrix}.
\end{equation}

(c) (Non-trapping condition) The initial data $u_0$ satisfies that the
bicharacteristic flow associated to \eqref{symbol} is non-trapping.

(d) (Growth condition of the first order coefficients) 
There exist $c_0>0$ and $N,\, \tilde M\in\zz^+$ large enough such that
\begin{equation}
\label{dddeca}
|\partial_x^{\alpha}b_{j_k}(x)|\leq \frac{c_0}{\langle x\rangle^N},\;\;\;|\alpha|\leq 
\tilde  M,\;\;\,j=1,2,\;\;k=1,..,n.
\end{equation}

(e) {Regularity of the coefficients} For $J\in \zz^+$ sufficiently large
$$
a_{jk},\,b_{m_j},\,c_m\in C^J_b(\rr^n:\cc),\;\;\;j,k=1,..,n, \;m=1,2.
$$
the following result was established in \cite{KePoRoVe}.

\begin{theorem}\cite{KePoRoVe}
 There exist $s,\,N>0$, $s>N$, and $\tilde N>0$, depending only on $n$, such that for $u_0\in H^s(\rr^n)\cap L^2(\langle 
x\rangle^Ndx)$
there exists $T=T(\|u_0\|_{s,2},\|\langle x\rangle ^{N/2}u_0\|_{2})>0$ such that 
the IVP \eqref{1.2semi-l}
has a unique solution $u$ defined in the time interval
$[0,T]$ satisfying
$$
u\in C([0,T]:H^s(\rr^n)\cap L^2(\langle x\rangle^Ndx))\equiv X_T^{s,N}.
$$
Moreover, 
$$
\int_0^T\int |J^{s+1/2}u(x,t)|^2\langle x\rangle^{-\tilde N}dxdt<\infty,
$$
and for every $ u_0\in H^s(\rr^n)\cap L^2(\langle x\rangle^Ndx)$ there exists a 
neighborhood $\mathcal U$ of $u_0$ and a
$T^{\prime}>0\,$ such that the map data $\to$ solution of is continuous from $\mathcal U\,$ into
$\,X^{s,N}_{T^{\prime}}$. 
\end{theorem}

Here  $\|\cdot\|_{s,2}$ denotes the norm in the Sobolev spaces 
$H^s(\rr^n)$.

 As in the previous semi-linear cases
the proof in \cite{KePoRoVe} was based on the contraction principle. 

 In this non-degenerate case the operator 
describing the \lq\lq integrating factor" in \eqref{intfac} 
has not been shown to be an $L^2$-bounded operator.
However, thanks to the local smoothing effects
it suffices to solve \eqref{commutator} up to \lq\lq small" first order term. This is achieved in 
\cite{KePoRoVe} by introducing a new class of symbols.

\begin{definition}
\label{Definition 1.2}
(i) It will be said that $a\in \mathcal S(\rr^n:S^m_{1,0})$ 
(where $S^m_{1,0}$ the class of $\psi$.d.o's of classical symbols 
of order $m$) if
$a\in C^{\infty}(\rr^n\times\rr^n\times \rr^n)$ 
and $a$ satisfies
\begin{equation}
\label{scw}
|\langle z\rangle^{\mu}\partial_s^{\alpha}\partial_x^{\beta}\partial_{\xi}^{\gamma}\,a(s;x,\xi)|
\leq c_{\mu \alpha \beta\gamma}\langle\xi\rangle^{m-|\gamma|},\;\;\forall\,z,x,\xi\in(\zz^+)^n,
\;\forall \,\mu,\alpha, \beta, 
\gamma\in(\zz^+)^n.
\end{equation}
(ii) For $a\in \mathcal S(\rr^n:S^m_{1,0})$ let 
\begin{equation}
\label{class}
b(x,\xi)=\chi(|\xi|)a(P(x,A_h\xi);x,\xi),
\end{equation}
where $P(y,z) = y-(y\cdot z)z/|z|^2$ for $y,z\in\rr^n, \;z\ne 0$, is
the projection of $y$ onto the hyperplane perpendicular to $z$, $A_h$ as in \eqref{matrix}, and $\chi\in 
C^{\infty}(\rr^n)$
with $\chi(l)=0$ for $|l|<1$ and $\chi(l)=1$ for $|l|>2$.
\end{definition}

 In fact, we showed in \cite{KePoRoVe} that it  suffices to have \eqref{scw} for sufficiently large
$\mu,\,\alpha,\,\beta,\,\gamma \in(\zz)^n$.

 We observe that if  $a\in \mathcal S(\rr^n:S^m_{1,0})$, then
$$
\partial^{\alpha}_{\xi}(\xi^{\beta}a(\cdot))\in\mathcal S(\rr^n:S^k_{1,0}),\;\;k=m+|\beta|-|\alpha|,
$$
i.e. this class is closed with respect to differentiation and multiplication in the $\xi$-variable.
We shall show below that this class is also closed with respect to differentiation in the $x$-variable.

 Roughly speaking, for  $r\in\zz^+$ large enough 
$$
\langle x\rangle^{-r} \chi(|\xi|)a(P(x,A_h\xi);x,\xi)=
\langle x\rangle^{-r} b(x,\xi), 
$$
is a symbol in the class $S^m_{1,0}$.

 In \cite{KePoRoVe}, we deduce several properties of operators with symbol in our class. 
These include their continuity from $H^m(\rr^n)$ to $L^2$ and  their composition rules with
classical differential operators $P(x,\partial_x)$ with decaying coefficients, i.e. with
$P(x,\partial_x)=\phi_{\alpha}(x)\partial_x^{\alpha}$ with $|x|^l\phi_{\alpha}$ bounded for $l\in\zz$ 
sufficiently large.

 The proof of the nonlinear results
concerning the IVP \eqref{1.2semi-l} relies on two key
linear estimates.
The first one is concerned with the
smoothing effect described for  solutions of the IVP \eqref{1.2linear} with
$(a_{jk}(x))$ being just an invertible matrix 
\begin{equation}
\label{selh}
\aligned
&\int_0^{T}\int_{\rr^n}|J^{s+1/2}u(x,t)|^2\langle x\rangle^{-\tilde N}dxdt
\\
&\leq c(1+T)\sup_{0\leq t\leq 
T}\|u(t)\|^2_{s,2}+c\int_0^{T}\int_{\rr^n}|f(x,t)|^2dx\,dt,
\endaligned
\end{equation}
or
\begin{equation}
\label{seli}
\begin{aligned}
&\int_0^{T}\int_{\rr^n}|J^{s+1/2}u(x,t)|^2\langle x\rangle^{-\tilde N}dxdt\\
&\leq c(1+T)\sup_{0\leq t\leq T}\|u(t)\|^2_{s,2} +
\int_0^{T}\int_{\rr^n}|J^{s-1/2}f(x,t)|^2\langle x\rangle^{\tilde N}dxdt
\end{aligned}
\end{equation}
for $\tilde N>1$.

The second is related with the local wellposedness in $L^2$ (and in $H^s$) of
the IVP \eqref{1.2linear}.
To establish it we  followed an indirect approach.
First, we truncated at infinity the
operator $\mathcal L(x)=- \partial_{x_k}a_{jk}(x)\partial_{x_j}$ using $\theta\in C^{\infty}_0(\rr^n)$
with $\theta(x)=1,\;|x|\leq 1$, and $\theta(x)=0,\; |x|\geq2$. For $R>0$ we define
\begin{equation}
\label{trunc}
\mathcal L^R(x)=\theta(x/R)\mathcal L(x)+(1-\theta(x/R))\mathcal L^0,
\end{equation}
where $\mathcal L^0=-a^0_{jk}\partial^2_{x_jx_k}$, $A_h=(a^0_{jk})_{j,k=1,..,n}$ 
as in \eqref{matrix}, with the
decay assumption $a_{jk}(x)-a_{jk}^0\in\mathcal S(\rr^n)$, $j,k=1,..,n$,
(the same proof worked if
one just assumed that
the corresponding estimate held for a sufficiently large number of semi-norms of $\mathcal S(\rr^n)$).
Thus,
$$
\mathcal L(x)=\mathcal L^R(x)+\mathcal E^R(x).
$$
For $R$ large enough we considered the bicharacteristic flow
$(X^R(s;x,\xi),\Xi^R(s;x,\xi))$ associated to the operator $\mathcal L^R(x)$ and  the
corresponding integrating factor $K^R$, i.e. the operator with symbol as in \eqref{intfac}
but evaluated in the
bicharacteristic flow $(X^R(s;x,\xi),\Xi^R(s;x,\xi))$.
We deduced  several estimates concerning the operator $K^R$.
In particular, we showed in \cite{KePoRoVe} that there exists $N_0\in\zz^+$ (depending only on the dimension $n$)
such that for any $M\in\zz^+$ there exist $N_1=N_1(M)\in\zz^+$ and $R_0$ sufficiently large
such that for $R\geq R_0$ it follows that  
\begin{equation}
\label{trunc1}
\begin{aligned}
&\sup_{0\leq t\leq T}\|K^R u(t)\|^2_2 \leq c_0 R^{N_0}\|u(0)\|_2^2\\
&+R^{-M}\int_0^T\int_{\rr^n}|J^{1/2}u|^2\langle x\rangle^{-\tilde N}dxdt+
c_0T R^{N_0+N_1(M)} \sup_{0\leq t\leq T}\|u(t)\|^2_2,
\end{aligned}
\end{equation}
for $\tilde N$ large enough depending only on the dimension $n$.
 In \cite{KePoRoVe} the coefficients in \eqref{1.2semi-l} were taken in the  
Schwartz class $\mathcal S(\rr^n)$. However, it is clear from the proof there that the same argument works 
if one just assumes that a fixed large number of semi-norms of $\mathcal S(\rr^n)$ of the
coefficients are bounded.
In this case, $M$ will be chosen depending only on the decay of the coefficients.
More precisely, one chooses $M=M(N)$, $M(N)\uparrow \infty$ as $N \uparrow \infty$, with $N$ as in  
\eqref{ddeca}, \eqref{dddeca}.

To complete  the estimate one needs to  consider the operator $E^R=I-\tilde K^R
(K^R)^*$,where the symbol of $\tilde K^R$ differs from that of $K^R$ only
in the sign of the exponent, and
$(K^R)^*$ is the adjoint of $K^R$. It was established   that $E^R u(t)$ satisfies
an estimate similar to that in \eqref{selh}.
Combining these results  we get that
\begin{equation}
\label{key1}
\begin{aligned}
&\sup_{0\leq t\leq T}\| u(t)\|^2_2 \leq c_0 R^{N_0}\|u(0)\|_2^2\\
&+R^{-M}\int_0^T\int_{\rr^n}|J^{1/2}u|^2\langle x\rangle^{-\tilde N}dxdt+
c_0T R^{N_0+N_1(M)} \sup_{0\leq t\leq T}\|u(t)\|^2_2\\
&+c_0R^{N_0}\int_0^T\int_{\rr^n}|f(x,t)|^2dxdt.
\end{aligned}
\end{equation}

 From \eqref{selh} and  \eqref{key1} fixing $T$ small enough one gets that
\begin{equation}
\label{key}
\sup_{0\leq t\leq T}\| u(t)\|^2_2
\leq c_0 R^{N_0}\|u(0)\|_2^2+c_0R^{N_0}\int_0^T\int_{\rr^n}|f(x,t)|^2dxdt.
\end{equation}

 This allows to use the contraction principle to obtain in \cite{KePoRoVe} the desired nonlinear result.

 Returning to the IVP \eqref{1.1} we shall assume that the
coefficients satisfy the following hypotheses: 

 (H1) \underbar{Non-degeneracy.} Given $r_0>0$ there exists $\gamma_{r_0}\in (0,1)$ such that for any 
\begin{equation} 
\label{H1}
\begin{aligned} 
&(x,t,\vec z)\in\rr^n\times \rr\times (\overline{(B_{r_0}(0))})^{2n+2}\equiv D_{r_0},\\ 
&\gamma_{r_0} |\xi|\leq |a_{jk}(x,t,\vec z)\xi|\leq \gamma_{r_0}^{-1}|\xi|,\;\;\;\;\forall \xi\in\rr^n, 
\end{aligned} 
\end{equation} 
where
$\overline{B_{r_0}(0)}=\{z\in\cc\,:\,|z|\leq r_0\}$, and with 
\begin{equation} 
\label{h1a} 
\begin{aligned} 
&A(x,0,\vec 0)=(a_{jk}(x,0,\vec 0))_{j,k=1,..,n}\\ 
&=A_0(x)+A_h=(a_{0,jk}(x))_{j,k=1,..,n}+(a^0_{jk})_{j,k=1,..,n} 
\end{aligned}
\end{equation} 
where for some $N,\, \tilde M\in \zz^+$ large enough 
\begin{equation} 
\label{H1b}
|\partial_x^{\alpha}a_{0,jk}(x)|\leq \frac{c_0}{\langle x \rangle^N},\;\;\;\;\;|\alpha|\leq 
\tilde M\;\;\;j,k=1,..,n,
\end{equation} 
and $ A_h$ as in \eqref{matrix}.

(H2)  \underbar{Asymptotic flatness.} There exists $c>0$ such that for any $ (x,t)\in \rr^n\times
\rr$ and $0<|\alpha|\leq \tilde M$, $\,0\leq|\alpha'|\leq \tilde M$ it follows that \begin{equation} \label{H2}
|\partial_{x}^{\alpha}a_{jk}(x,t,\vec 0)| + |\partial_t\partial_{x}^{\alpha'}a_{jk}(x,t,\vec 0)|\leq \frac{c}{\langle
x\rangle^N}, \end{equation} for $k,j=1,..,n$ with $N,\,\tilde M $ as in \eqref{H1b}.

(H3) \underbar {Growth of the first order coefficients.}  There exist $c,\,c_0>0$ such that
for any $x\in\rr^n$ and any $\,(x,t)\in \rr^n\times \rr$
\begin{equation}
\label{H3}
\begin{cases}
\begin{aligned}
&| \partial_{x}^{\alpha} b_{m_j}(x,0,\vec 0)|\leq \frac{c_0}{\langle 
x\rangle^N},\;\;|\alpha|\leq \tilde 
M,\;\;m=1,2,\;j=1,..,n,\\
&|\partial_t\partial_{x}^{\alpha}  b_{m_j}(x,t,\vec 0)|\leq \frac{c}{\langle x\rangle^N}
,\;\;|\alpha|\leq \tilde
M,\;\;m=1,2,\;j=1,..,n,
\end{aligned}
\end{cases}
\end{equation}
where $\vec b_l=(b_{l_1},..,b_{l_n}),\,l=1,2$  with $N,\,\tilde M $ as in \eqref{H1b}.

(H4) \underbar {Regularity.} For any $J\in\zz^+$ and  $ r_0>0$ the coefficients
\begin{equation}
\label{H4a}
a_{jk},\,b_{1_j},\,b_{2_j}
\in C_b^J(\rr^n\times \rr\times (\overline{(B_{r_0}(0))})^{2n+2}),
\end{equation}
and
\begin{equation}
\label{H4b}
c_1,\,c_2\in C_b^J(\rr^n\times \rr\times (\overline{(B_{r_0}(0))})^{2}).
\end{equation}

 Our main result is the following.

\begin{theorem}
\label{Theorem A}

 Under the hypotheses (H1)-(H4) there exists $N=N(n)\in \zz^+$ such that given any
\begin{equation}
\label{clase}
u_0\in H^s(\rr^n)\;\;\;\;\text{with}\;\;\;\;\langle 
x\rangle^N\partial_x^{\alpha}u_0\in L^2(\rr^n),
\;\,|\alpha|\leq s_1,
\end{equation}
and
\begin{equation}
\label{clase2}
f\in L^{1}(\rr^+: H^s(\rr^n))\,\;\;\;\;\;\text{with}\;\;\;\;\;\langle 
x\rangle^N\partial_x^{\alpha}f\in L^{1}(\rr^+: L^2(\rr^n)),
\;\,|\alpha|\leq s_1,
\end{equation}
where $s, s_1\in\zz^+$, sufficiently large with $s>s_1+4$, for which the Hamiltonian flow
$H_{h(u_0)}$ associated to the symbol
\begin{equation}
\label{H5}
h(u_0)=
h_{u_0}(x,\xi)=\sum_{j,k=1}^na_{jk}(x,0,u_0,\bar u_0,\nabla 
u_0,\nabla \bar
u_0)\xi_j\xi_k,
\end{equation}
is non-trapping, there exist $T_0>0$, depending on
$$
\aligned
\lambda&=\|u_0\|_{s,2}+\sum_{|\alpha_1|\leq s_1}\|\langle 
x\rangle^N\partial_x^{\alpha}u_0\|_2\\
&+\int_0^{\infty}\,\|f(t)\|_{s,2}\,dt
+\sum_{|\alpha_1|\leq s_1}\,\int_0^{\infty}\,\|\langle 
x\rangle^N\partial_x^{\alpha}f(t)\|_2dt,
\endaligned
$$
the constants in (H1)-(H4) and on the non-trapping condition 
(H5), and a unique solution
$u=u(x,t)$ of the equation in \eqref{1.1} with initial data $u(x,0)=u_0(x)$ on the time
interval
$[0,T_0]$ satisfying

\begin{equation}
\label{1.27a}
\begin{aligned}
&u\in C([0,T_0]:H^{s-1})\cap 
L^{\infty}([0,T_0]:H^{s})\cap   
C^1((0,T_0):H^{s-3}),\\
&\langle
x\rangle^N\partial_x^{\alpha}u\in 
C([0,T_0]:L^2),\;\;\,|\alpha|\leq 
s_1.
\end{aligned}
\end{equation}
Moreover, if $(u_0,f)\in 
H^{s'}(\rr^n)\times L^{1}(\rr^+ : H^{s'}(\rr^n))$
with $s'>s$ then \eqref{1.27a} holds with $s'$ instead of $s$ in the same time interval 
$[0,T_0]$.

\end{theorem}

\underbar{Remarks}   Here, we are not concerned with the problem of estimating 
the optimal values of $s,\,s_1$ or $N$ in Theorem \ref{Theorem A}. 

Also
we shall not attempt to obtain the sharp form
of the persistence property of the solution,
(i.e. $u\in C([0,T_0]:H^{s}(\rr^n))$) as well as the 
continuous dependence of the solution upon the data. These can be established by combining
the argument in \cite{BoSm} with our key \it a  priori \rm estimates in Lemma 2.1.

 Similarly, from our arguments it is easy to deduce that the local solution possesses the local 
smoothing effect, i.e. if $u_0\in  H^{s}(\rr^n)$, then $J^{s+1/2}u\in
L^2(\rr^n\times [0,T_0]:\langle x\rangle^{-\tilde N}dxdt)$.

 The use of the weights in  Theorem \ref{Theorem A} comes from two sources.
First, in order to convert $L^1$ conditions, such as \eqref{ichi}, into $L^2$ conditions.
Secondly, one needs them in order to maintain the asymptotic flatness condition (H2), when 
$\vec 0$ is replaced 
by $(u,\bar u,\nabla_x u,\nabla_x \bar u)$, for a solution $u$.

\vskip.01in

As it was mentioned quasilinear results as those in Theorem \ref{Theorem A}
cannot be obtained by using just a fixed point argument. 
Instead, as in \cite{KePoVe4},  we shall rely on the artificial viscosity method.
First, we consider the linear problem
\begin{equation}
\begin{cases}
\label{linearvis}
\begin{aligned}
&\partial_t u=-\epsilon \Delta^2u-i\,\partial_{x_j}(a_{jk}(x,t) \partial_{x_k}u)+\vec b_1(x,t)\cdot \nabla u
+\vec b_2(x,t)\cdot \nabla \bar u\\
&+ c_1(x,t)u + c_2(x,t) \bar u +f(x,t)\equiv -\epsilon \Delta^2u + L(x,t)u+f(x,t),\\
&u(x,0)=u_0(x),
\end{aligned}
\end{cases}
\end{equation}
with $\epsilon\in(0,1)$.
 The main step is to obtain the following \it a priori \rm estimate
for solutions of the linear IVP \eqref{linearvis} : there exits $T>0$ such that  
\begin{equation}
\label{estimate}
\sup_{0\leq t\leq T}\| u(t)\|^2_2\leq c_T(\|u_0\|^2_2+\int_0^T 
\|f(t)\|^2_2dt),
\end{equation}
with $T,\,c$ independent of $\epsilon$
and depending on an appropriate manner on the coefficients in \eqref{linearvis}. 

 The inequality \eqref{estimate} will be proved under general hypotheses
on the coefficients in \eqref{linearvis}. This will allow us to find a class of functions such that
when the coefficients of the equation in \eqref{1.1} are evaluated in an element in this class
they satisfy these general hypotheses.

 We denote by $u^{\epsilon}$ the solution of the following nonlinear IVP 
associated to equation in \eqref{1.1}
\begin{equation}
\begin{cases}
\label{111}
\begin{aligned}
&\partial_t u=-\epsilon \Delta^2u-\, i\,\partial_{x_j}(a_{jk}(x,t,u,\bar u,\nabla u,\nabla \bar u) \partial_{x_{jk}}u)\\
&+\vec b_1(x,t,u,\bar u, \nabla u,\nabla \bar u)\cdot \nabla u
+\vec b_2(x,t,u,\bar u, \nabla u,\nabla \bar u)\cdot \nabla \bar u\\
&+ c_1(x,t, u,\bar u)u + c_2(x,t, u,\bar u) \bar u +f(x,t),\\
&u(x,0)=u_0(x).
\end{aligned}
\end{cases}
\end{equation}

 The viscosity method provides the solution $u^{\epsilon}$, in the appropriate class, 
in the time interval $[0,T_{\epsilon}]$ with $T_{\epsilon}=O(\epsilon)$. Evaluating 
the coefficients in \eqref{111} in this solution we get a linear problem as that in 
\eqref{linearvis} for which 
the estimate \eqref{estimate} holds. This  \it a priori \rm estimate allows us
to extend the solution $u^{\epsilon}$, in the same class, to a time interval
$[0,T_0]$ with $T_0$ independent of $\epsilon\in (0,1)$. 

  Once the estimate \eqref{estimate} is available the proof of Theorem \ref{Theorem A} follows 
an argument quite similar to that explained 
in detail in \cite{KePoVe4}. So we shall concentrate in the proof 
of the inequality  \eqref{1.1}. This will be given in section 2, Lemma \ref{mainlth}.

Finally, we point out another difference between the elliptic quasi-linear case and 
the non-degenerate one.

 In \cite{KePoVe4} the following general class of quasilinear equation was considered
\begin{equation}
\begin{cases}
\label{geneq}
\begin{aligned}
&\partial_t u= - i\,\partial_{x_j}(a_{jk}(x,t,u,\bar u,\nabla u,\nabla \bar u) \partial_{x_k}u
+\partial_{x_j}(b_{jk}(x,t,u,\bar u,\nabla u,\nabla \bar u) \partial_{x_k}u\\
&+\vec b_1(x,t,u,\bar u, \nabla u,\nabla \bar u)\cdot \nabla u
+\vec b_2(x,t,u,\bar u, \nabla u,\nabla \bar u)\cdot \nabla \bar u\\
&+ c_1(x,t, u,\bar u)u + c_2(x,t, u,\bar u) \bar u +f(x,t),
\end{aligned}
\end{cases}
\end{equation} 
where $(b_{jk})$ is a symmetric complex valued matrix.

 Under the ellipticity assumption : given $r_0>0$ there exists $\gamma_{r_0}>0$ such that for any
$$
(x,t,\vec z)\in\rr^n\times \rr\times \overline{((B_{r_0}(0))}^{2n+2},
$$
it follows that 
$$
\langle a_{jk}(x,t,\vec z)\xi,\xi\rangle-|\langle b_{jk}(x,t,\vec z)\xi,\xi\rangle|
\geq \gamma_{r_0} |\xi|^2,\;\;\forall \xi\in\rr^n,
$$
with $B_{r_0}(0)$ as in \eqref{H1} and  similar assumptions
on the asymptotic flatness,
the growth of first order coefficients, the  regularity of the coefficients, and a non-trapping
hypothesis it was established in \cite{KePoVe4} that the IVP associated to the equation 
\eqref{geneq} is locally well posed.

The non-trapping
hypothesis in \cite{KePoVe4} was the following :
Given $u_0\in H^r(\rr^n), r>n/2+2$, define
$$
\varpi_2(x,0,\xi)=-a_{jk}(x,0,u_0,\bar u_0,\nabla u_0,\nabla \bar u_0)\xi_k\xi_j,
$$
$$
\varpi_3(x,0,\xi)=-b_{jk}(x,0,u_0,\bar u_0,\nabla u_0,\nabla \bar u_0)\xi_k\xi_j,
$$
and
$$
\kappa(x,\xi)=\sqrt{\varpi^2_2(x,0,\xi)-|\varpi_3(x,0,\xi)|^2}.
$$
 It was said that $u_0$ satisfies the non-trapping hypothesis  if there exists
$0<\eta<1$ and functions $a(x, \xi), a_1(x, \xi)$ such that
$$
\kappa(x,\xi)=a(x,\xi)+\eta a_1(x,\xi),
$$
(with $a(x,\xi)$ real, homogeneous of degree $2$, with $|\partial_x^{\beta}a(x,\xi)|\in C^{1,1}(\rr^n\times\rr^n)$,
 $|\beta|\leq N(n)$  and $\theta(\xi)\,a(x,\xi)\in C^{N(n)}(\rr^n\times\rr^n)$, where
$\theta\equiv 1$ for $|\xi|>1$, $\theta\equiv 0$ for $|\xi|<1/2$, $\theta\in C^{\infty}(\rr^n)$, with $a_1$
verifying similar estimates), the Hamiltonian flow $H_a$ associated to the symbol $a$ is non-trapping, for 
more details 
see 
\cite{KePoVe4}.

 In the non-degenerate setting the equation in \eqref{geneq} does not allow for a more general one that that  considered
in \eqref{1.1}. This is due to the following linear algebra result (whose proof follows by induction).

\begin{lemma}
\label{Lemma 1.3}

Let $A,\;B$ be two $n\times n$ matrices, $A$ be real symmetric non-positive and non-degenerate, and
$B$ a symmetric complex valued one such that
$$
\langle A\xi,\xi\rangle=0\;\;\;\text{implies}\;\;\;
\langle B\xi,\xi\rangle=0.
$$
Then there exists $\lambda \in\cc$ such that $B=\lambda A$.
\end{lemma}

\section{The linear problem}
\label{motion}

 In this section we shall consider the linear IVP
\begin{equation}
\begin{cases}
\label{1.2linear2}
\begin{aligned}
&\partial_t u=- i\,\partial_{x_j}(a_{jk}(x,t) \partial_{x_k}u)+\vec b_1(x,t)\cdot 
\nabla u
+\vec b_2(x,t)\cdot \nabla \bar u\\
&+ c_1(x,t)u + c_2(x,t) \bar u +f(x,t)\equiv L(x,t)u+f(x,t),\\
&u(x,0)=u_0(x),
\end{aligned}
\end{cases}
\end{equation}
where $x\in\rr^n,\,n>1$, $t\in[0,T]$, with $T>0$, and its associated $\epsilon$-viscosity version
\begin{equation}
\label{1.3}
\begin{cases}
\begin{aligned}
&\partial_t u=-\epsilon \Delta^2u+ L(x,t)u+f(x,t),\;\;\;\;\epsilon\in(0,1),\\
&u(x,0)=u_0(x),
\end{aligned}
\end{cases}
\end{equation}
under the following assumptions:

(H$_l$1) \underbar{Non-degeneracy} : $A(x,t)=(a_{jk}(x,t))_{j,k=1,..,n}$ is a real 
symmetric matrix and there exist $ \gamma,\,\gamma_0\in(0,1)$ such that for any 
$\xi\in \rr^n, (x,t)\in 
\rr^n\times [0,T]$
\begin{equation}
\label{1.4}
\begin{cases}
\begin{aligned}
&\gamma|\xi|\leq |A(x,t)\xi|\leq \gamma^{-1}|\xi|,\\
\\
&\gamma_0|\xi|\leq |A(x,0)\xi|\leq \gamma_0^{-1}|\xi|,
\end{aligned}
\end{cases}
\end{equation}
with
\begin{equation}
\label{1.5}
\begin{aligned}
&A(x,0)=(a_{jk}(x,0))_{j,k=1,..,n}\\
&=A_0(x)+A_h=(a_{0,jk}(x))_{j,k=1,..,n}+(a^0_{jk})_{j,k=1,..,n}
\end{aligned}
\end{equation}
where for some $N,\,\tilde M\in \zz^+$ large enough
\begin{equation}
\label{1.6}
|\partial_x^{\alpha}a_{0,jk}(x)|\leq \frac{c_0}{\langle x \rangle^N},\;\;\;\;\;|\alpha|\leq 
\tilde  M,\;\;\;j,k=1,..,n,
\end{equation}
and
$A_h$ as in \eqref{matrix}.

(H$_l$2) \underbar{Asymptotic Flatness} : There exists $ c>0$ such that for all $\xi\in \rr^n, (x,t)\in 
\rr^n\times [0,T]$, $j,k=1,..,n$, and $0<|\alpha|\leq \tilde M,\; 0\leq|\alpha'|\leq  \tilde M$
\begin{equation}
\label{1.8}
|\partial_x^{\alpha}a_{jk}(x,t)|+|\partial_t\partial_x^{\alpha'}a_{jk}(x,t)|\leq 
\frac{c}{\langle x 
\rangle^N},
\end{equation}
with $N,\, \tilde M$ as in $H_l1$.

(H$_l$3) \underbar{Growth Condition of the First Order Coefficients} : There exist $ c,\;c_0>0$ such that  
for all $\xi\in \rr^n, (x,t)\in
\rr^n\times [0,T]$
\begin{equation}
\label{1.9.a}
|\partial_x^{\alpha}b_{m_j}(x,0)| \leq \frac{c_0}{\langle x
\rangle^N},\;\;\;|\alpha|\leq \tilde M,\;\;m=1,2,\;\;j=1,..,n,
\end{equation}
and
\begin{equation}
\label{1.9.b}
|\partial_x^{\alpha}b_{m_j}(x,t)| \leq \frac{c}{\langle x
\rangle^N},\;\;\;|\alpha|\leq \tilde M,\;\;m=1,2,\;\;j=1,..,n,\\
\end{equation}
with $N,\, \tilde M$ as in $H_l1$.

(H$_l$4) \underbar{Regularity of the Coefficients} :
\begin{equation}
\label{1.10}
a_{jk},\,b_{m_j},\,c_m\in C^J_b(\rr^n\times [0,T]),\;\;\;j,k=1,..,n, \;m=1,2,
\end{equation}
with $J=J(n)\in \zz^+$ sufficiently large such that the proofs below involving $\psi$.d.o's
can be carried out, with norm
\begin{equation}
\label{1.11}
\|\partial_x^{\alpha}\partial_t^r d\|_{L^{\infty}(\rr^n\times [0,T])}\leq c_J,\;\;\;|\alpha|+r\leq J.
\end{equation}
with $d=a_{jk},\,b_{m_j}$ or $c_m$, $j,k=1,..,n, \;m=1,2$.

(H$_l$5) \underbar{Non-trapping Condition} : The bicharacteristic flow associated to 
the symbol of $A(x,0)$ (see \eqref{1.5}), i.e. 
\begin{equation}
\label{symsym}
h_2^0(x,\xi)=\sum_{j,k=1}^n a_{jk}(x,0)\xi_j\xi_k
\end{equation}
 is non-trapping.

 We shall use $c_0$ to denote a generic constant  which only depends on  the coefficients evaluated
at time $t=0$.

The main ingredient in the proof of our main result Theorem \ref{Theorem A} is the 
following estimate for the solution of the linear IVP
\eqref{1.2linear2}.

\begin{lemma}
\label{mainlth}
 There exist $\tilde c_0=\tilde c_0(\gamma_0,c_0,n)>0$ with $\gamma_0$ as in \eqref{1.4}, 
$c_0$ as in \eqref{1.6} and \eqref{1.9.a}, 
and 
 the constant $c_0$ in Lemma \ref{tidoi}  below, $\tilde N=\tilde N(n)$, and
$K_0>0, \;T_0\in(0,T]$ depending on  
\eqref{1.4}-\eqref{1.11}  
such that for $\tilde T\in (0,T_0)$ the solution of the IVP \eqref{1.3} $u^{\epsilon}$ satisfies
\begin{equation}
\label{1.12}
\begin{aligned}
&\sup_{0\leq t\leq \tilde T}\|u^{\epsilon}(t)\|_2+\int_0^{\tilde T}\int|J^{1/2}u^{\epsilon}|^2\langle 
x\rangle^{-\tilde N}
dxdt\\
&\leq \tilde c_0 e^{K_0\tilde T}(\|u_0\|_2+\int_0^{\tilde T}\|f(\cdot,t)\|_2dt).
\end{aligned}
\end{equation}

Moreover, \eqref{1.12} still holds if we replace its last term by
$$
(\int_0^{\tilde T}\int|J^{-1/2}f(x,t)|^2\langle
x\rangle^{-\tilde N} dxdt)^{1/2}.
$$
\end{lemma}

We recall the class of $\psi$.d.o's introduced in \cite{KePoRoVe}:
$$
\Psi_af(x)=\int e^{ix\cdot \xi}a(x,\xi)\hat f(\xi)d\xi,
$$
with symbol
$$
a(x,\xi)=\chi(|\xi|)a(P(x,A_h\xi);x,\xi),
$$
with $\chi$ as in \eqref{class}, 
$$
P(x,A_h\xi)=x-\bigl(\frac{x\cdot A_h\xi}{|\xi|^2}\bigr)A_h\xi,
$$
i.e. $P(x,A_h\xi)$ is the projection of $x$ into the hyperplane perpendicular to $A_h\xi$, and
$$
a=a(s;x,\xi)\in \mathcal {S}(\rr^n:S^m_{1,0}),
$$
$\mathcal {S}(\cdot)$ denoting the Schwartz class, and $S^m_{1,0}$  the class of classical symbol of $\psi$.d.o's
of order $m$.

 In \cite{KePoRoVe} (Theorem 3.2.1) we prove that
\begin{equation}
\label{1.13}
\|\Psi_a f\|_2\leq c \|f\|_{m,2}.
\end{equation}

We shall use that our class is closed under differentiation of the symbol in the
$x$-variable.

\begin{lemma}
\label{ccd}

Let
$$
a_{\alpha}(x,\xi)=\partial_x^{\alpha}(a(P(x,A_h\xi);x,\xi)\chi(|\xi|)),
$$
with $a(s;x,\xi)$ as above. Then $a_{\alpha}(x,\xi)$ defines a symbol in our class. Moreover,
\begin{equation}
\label{1.14}
\|\Psi_{a_{\alpha}} f\|_2\leq c_{\alpha} \|f\|_{m,2}.
\end{equation}
\end{lemma}

\begin{proof}
 
 First we consider the case $\alpha=(1,0,..,0)$. So
$$
\aligned
&\partial_{x_1}{a(P(x,A_h\xi);x,\xi)}\\
&=\frac{\partial a}{\partial x_1}(P(x,A_h\xi);x,\xi)
+\sum_{j=1}^n\frac{\partial a}{\partial z_j}(P(x,A_h\xi);x,\xi) \frac{\partial P}{\partial x_1}
(x,A_h\xi).
\endaligned
$$

Since
$$
\partial_{x_1}(P(x,A_h\xi))_j=\partial_{x_1}\bigl(x_j-\frac{x\cdot A_h\xi}{|\xi|^2}(A_h\xi)_j\bigr)=
\delta_{1j}-\frac{(A_h\xi)_1(A_h\xi)_j}{|\xi|^2},
$$
it follows that (see  remark 3.1.3(b) in \cite{KePoRoVe})
$$
\partial_{x_1}{a(P(x,A_h\xi);x,\xi)\chi(|\xi|)}=b_1(P(x,A_h\xi);x,\xi)\chi(|\xi|)),
$$
with $b_1=b_1(z;x,\xi)\in \mathcal {S}(\rr^n:S^m_{1,0})$ which yields the result.

The proof of the general case combines the above argument and induction in $|\alpha|$.
\end{proof} 

 Using  the notation
\begin{equation}
\label{1.15}
\mathcal{L}=\mathcal{L}(x,t)=-\partial_{x_j}(a_{jk}(x,t)\partial_{x_k}\cdot),
\end{equation}
and taking complex conjugate in the equation \eqref{1.3}
we obtain the system
\begin{equation}
\label{1.16}
\begin{cases}
\begin{aligned}
&\partial_t \vec w=-\epsilon \Delta^2 I\vec w + iH\vec w+B\vec w+C\vec w+\vec F,\\
&\vec w(x,0)=\vec w_0(x),
\end{aligned}
\end{cases}
\end{equation}
where
\begin{equation}
\label{1.17}
\vec w=\begin{pmatrix} u\\ \bar u\end{pmatrix},\;\;\;\; \vec F=\begin{pmatrix} f\\ \bar f\end{pmatrix},
\end{equation}
\begin{equation}
\label{1.18}
H=\begin{pmatrix} \mathcal{L}& 0\\0&-\mathcal{L}\end{pmatrix},\;\;\;\;
B=\begin{pmatrix} \Psi_{b_1}&\vec b_2(x,t)\cdot\nabla\\ \overline{\vec b_2(x,t)}\cdot\nabla &\Psi_{\bar 
b_1}\end{pmatrix},
\end{equation}
with $b_1\in S^{1}_{1,0}$ with odd symbol in $\xi$, and $C$ a $2\times 2$ matrix of $\psi$.d.o's of
order zero.

 At some point we will take derivative of the equation in \eqref{1.1}, so the new coefficients of $\nabla u$
will be a combination of the original ones and some derivatives of the $a_{jk}$'s. In fact,
these coefficients depend on the order of the derivative just as a multiplicative constant. For this reason we consider a 
general $b_1\in S^{1}_{1,0}$ with odd symbol in $\xi$ in \eqref{1.15}, \eqref{1.18}.

\begin{lemma}
\label{lse}

There exists $N\in\zz^+$ such that for any $T_0\in(0,T]$ and any $\epsilon\in (0,1)$
the solution $\vec w^{\epsilon}=\vec w$ of the IVP
\eqref{1.16} satisfies
\begin{equation}
\label{1.19}
\begin{aligned}
&\int_0^{T_0}\int |J^{1/2}\vec w|^2 \langle x\rangle^{-\tilde N}dxdt
\leq (c_0+cT_0)\sup_{0\leq t\leq T_0}\|\vec w(t)\|_2^2\\
&+c_0\epsilon \int_0^{T_0}\|\Delta \vec w(t)\|_2^2dt +c_0 \int_0^{T_0}\|\vec F(t)\|_2^2dt,
\end{aligned}
\end{equation}
where $c_0$ depends only on the coefficients evaluated at $t=0$  and on Lemma 2.1 and $c$ depends on the 
estimates in \eqref{1.4}-\eqref{1.11}.

\end{lemma}

To prove Lemma \ref{lse} we will follow the argument in \cite{KePoRoVe} (Lemma 5.1.3). First, we
recall Lemma 5.1.1 in \cite{KePoRoVe}. 
 
\begin{lemma}
\label{tidoi}
Let $A(x,0)$ be as in \eqref{1.5}. Assume  that the bicharacteristic flow is non-trapped, i.e.
$$
\{X(s;x_0,\xi_0)\,:\,s\in\rr^+\}
$$
is unbounded for any $(x_0,\xi_0)\in\rr^n\times\rr^n-\{0\}$.
Suppose 
$\lambda\in L^1([0,\infty))\cap C([0,\infty))$ is strictly positive and
non-increasing. Then there exist $c_0>0$ and a real symbol $p\in S^0_{1,0}$, both depending on 
$h^0_2(x,\xi)$ in \eqref{symsym} and $\lambda$
such that
$$
H_{h^0_2}p=\{h^0_2,p\}(x,\xi)\geq \lambda(|x|)\,|\xi|-c_0,\;\;\;\;\forall \,(x,\xi)\in
\rr^n\times\rr^n.
$$
\end{lemma}

 We fix $N$ as in \eqref{1.6}, \eqref{1.8}, \eqref{1.9.a}, and \eqref{1.9.b}, and choose 
$$
\lambda(\rho)=1/(1+\rho^2)^{\tilde N/2}=\langle \rho\rangle^{\tilde N}.
$$
to obtain the following time dependent version of Lemma \ref{tidoi}.

\begin{lemma}
\label{tddoi}

These exists $T_0>0$ depending only on H$_l$2 such that for any $(x,\xi)\in\rr^n\times \rr^n$
\begin{equation}
\label{1.20}
H_{h_2}p=\{h_2;p\}(x,t,\xi)\geq \frac{|\xi|}{2\langle x\rangle^{\tilde N}}-2c_0,
\end{equation}
where
\begin{equation}
\label{1.21}
h_2=h_2(x,t,\xi)=a_{jk}(x,t)\xi_k\xi_j,
\end{equation}      
and $p\in S^0_{1,0}$ is the operator in Lemma \ref{tidoi}.
\end{lemma}

\begin{proof}

 We  recall the notation
$$
h^0_2(x,\xi)=a_{jk}(x,0)\xi_k\xi_j.
$$
 
By definition 
$$
\aligned
H_{h_2}p&=\partial_{\xi_j}h_2\partial_{x_j}p-\partial_{x_j}h_2\partial_{\xi_j}p\\
&=2a_{jk}(x,t)\xi_k\partial_{x_j}p-\partial_{x_j}a_{jk}(x,t)\xi_k\xi_j\partial_{\xi_j}p.
\endaligned
$$

From our hypothesis H$_l$2 it follows that for any $t\in(0,T_0]$ 
$$
\aligned
&|\partial_{\xi_j}h_2(x,t,\xi)-\partial_{\xi_j}h^0_2(x,\xi)|\\
&=2|(a_{jk}(x,t)-a_{jk}(x,0))\xi_k|\leq
c \frac{T_0}{\langle x\rangle^N} |\xi|,
\endaligned
$$
and
$$
\aligned
&|\partial_{x_j}h_2(x,t,\xi)-\partial_{x_j}h^0_2(x,\xi)|\\
&=|(\partial_{x_j}a_{jk}(x,t)-\partial_{x_j}a_{jk}(x,0))\xi_k\xi_j|
\leq
c \frac{T_0}{\langle x\rangle^N} |\xi|^2.
\endaligned
$$

Since $p\in S^0_{1,0}$ choosing $N\gg \tilde N$ one has that
$$
|H_{h_2}p-H_{h_2^0}p|\leq c \frac{T_0}{\langle x\rangle^{\tilde N}} |\xi|.
$$
 
From Lemma \ref{tidoi} we know that
$$
H_{h_2^0}p(x,\xi)\geq \frac{|\xi|}{\langle x\rangle^{\tilde N}} -c_0,
$$
so taking $T_0$ sufficiently small we obtain the desired result \eqref{1.20}.
\end{proof}

Now we shall prove Lemma \ref{lse}. 

\begin{proof} 

Let
\begin{equation}
\label{1.22}
k(x,\xi)=\begin{pmatrix}
exp(p(x,\xi))&0\\0&-exp(p(x,\xi))
\end{pmatrix},
\end{equation}
where $p\in S^0_{1,0}$ is the symbol in Lemma \ref{tddoi} 
so $K=\Psi_k$ is a diagonal matrix of $\psi$.d.o's of
order zero. We calculate
\begin{equation}
\label{1.23}
\begin{aligned}
&\partial_t \langle K\vec w,\vec w\rangle = \langle K\partial_t\vec w,\vec w\rangle + \langle K\vec w,\partial_t\vec w\rangle\\
&=-\epsilon(\langle K \Delta^2I\vec w,\vec w\rangle
+\langle K\vec w,\Delta^2I\vec w\rangle)\\
&+\langle (i[KH-HK]+KB+B^*K)\vec w,\vec w\rangle\\
&\langle (KC+C^*K)\vec w,\vec w\rangle
+\langle K\vec F,\vec w\rangle
+\langle K\vec w,\vec F\rangle.
\end{aligned}
\end{equation}

 We disregard the symbols of order zero and consider first the symbol of $i[KH-HK]+KB+B^*K$, i.e.
$$
\sigma(i[KH-HK]+KB+B^*K),
$$
which is equal up to a symbol of order zero to
\begin{equation}
\label{1.24}
-e^{p}\begin{pmatrix}
H_{h_2}p&0\\
0&H_{h_2}p\end{pmatrix}
+e^{p}\begin{pmatrix}
b_1(x,t,\xi)&2i\vec b_2(x,t)\cdot \xi\\
-2i\vec \bar b_2(x,t)\cdot \xi
&\overline{b_1(x,t,\xi)}\end{pmatrix}
\end{equation}
 We write
$$
\vec b_2(x,t)=\vec b_2(x,0)+t\;\frac{\vec b_2(x,t)-\vec b_2(x,0)}{t},
$$
with a similar identity for $  b_1(x,t,\xi)$. Also since
$K=\Psi_k$ has order zero it is easy to see that
$$
|-\epsilon(\langle K \Delta^2I\vec w,\vec w\rangle
+\langle K\vec w,\Delta^2I\vec w\rangle)|\leq c_0\epsilon \|\vec w(t)\|_{2,2}^2.
$$

 This combined with the matrix version of the sharp Garding inequality,
the hypothesis H$_l$3, and Lemma \ref{tddoi} yields, after integration by part, the desired result, i.e.
inequality \eqref{1.19} (for details we refer to the proof  of Lemma 5.1.3 
in \cite{KePoRoVe}).

\end{proof}

 Next, we shall recall some notations and definitions used in \cite{KePoRoVe}. First, we define
\begin{equation}
\label{truncado1}
a_{jk}^R(x,0)=\theta\bigl(\frac{x}{R}\bigr)a_{jk}(x,0)+\bigl(1-\theta\bigl(\frac{x}{R}\bigr)\bigr)a^0_{jk},
\end{equation}
where $\theta\in C^{\infty}_0(\rr^n)$, $\theta(x)=1$ for $|x|<1$, and $\theta(x)=0$ for $|x|>2$,
\begin{equation}
\label{truncado}
\mathcal{L}^R(x)=\partial_{x_j}(a_{jk}^R(x,0)\partial_{x_k}),\;\;\text{and}\;\;
\mathcal{L}(x,t)=\partial_{x_j}(a_{jk}(x,t)\partial_{x_k}).
\end{equation}

 We apply the operator $(K^R)^*$, which is independent of $t$, to the equation in \eqref{1.3} and use that
$$
\mathcal{L}(x,t)=\mathcal{L}(x,t)-\mathcal{L}(x,0)+\mathcal{L}^R(x)+\mathcal{E}^R(x),
$$
with
$$
\mathcal{E}^R(x)=\partial_{x_j}\bigl(\bigl(1-\theta\bigl(\frac{x}{R}\bigr)\bigr)(a_{jk}(x,0)-a^0_{jk})\partial_{x_k}\bigr),
$$
and
$$
\vec b_j(x,t)\cdot\nabla=(\vec b_j(x,t)-\vec b_j(x,0))\cdot\nabla+\vec b_j(x,0)\cdot\nabla, \;\;\;j=1,2.
$$
 
 We shall need the following symbols used in \cite{KePoRoVe},
\begin{equation}
\label{list}
\begin{aligned}
&b_1(x,\xi)=-b_1(x,0)i\xi,\\
&p^R(x,\xi)=- \chi(|\xi|)\int_{-\infty}^0 b_1(X^R(s;x,\xi),\Xi^R(s;x,\xi))ds,
\\
&p_e^R(x,\xi)=\frac 12(p^R(x,\xi)+p^R(x,-\xi)),
\\
&k^R(x,\xi)=\exp (p_e^R(x,\xi)),
\\
&k^R(x,\xi)=\exp(-p_e^R(x,\xi)).
\end{aligned}
\end{equation}
where 
$$
(X^R(s;x,\xi),\Xi^R(s;x,\xi)),
$$
denotes the bicharacteristic flow associated to the symbol of the truncated operator $\mathcal{L}^R(x)$
defined  in \eqref{truncado} and $\chi$ as in \eqref{class}.

We remark that after taking $s$-derivatives in the equation \eqref{1.1}  
and linearizing the resulting equation the new term $b_1(x,\xi)$ obtained has the form
$$
s\,\sum\limits_{j,k,l}\partial_{x_j}a_{jk}^R(x,0)
\xi_j\xi_k\xi_l\langle\xi\rangle^{-2}-\,b_1(x,0)i \xi.
$$

From our hypotheses it is clear that this new $b^R(x,\xi)$ satisfies similar  
estimates than that in
\eqref{list}.

Thus, from \eqref{1.3} we get
\begin{equation}
\label{1.25}
\begin{aligned}
&\partial_t(K^R)^*u=-\epsilon (K^R)^* \Delta^2u - i[(K^R)^*;\mathcal{L}^R]u+i\mathcal{L}^R(K^R)^*u\\
&+i[(K^R)^*;(\mathcal{L}(x,t)-\mathcal{L}(x,0))]u+i(\mathcal{L}(x,t)-\mathcal{L}(x,0))(K^R)^*u\\
&+i\mathcal{E}^R(K^R)^*u+i[\mathcal{E}^R;(K^R)^*]u\\
&+(K^R)^*\vec b_1(x,0)\cdot\nabla u + (K^R)^*(\vec b_1(x,t)-\vec b_1(x,0))\cdot\nabla u\\
&+(K^R)^*\vec b_2(x,0)\cdot\nabla \bar u + (K^R)^*(\vec b_2(x,t)-\vec b_2(x,0))\cdot\nabla \bar u\\
&+(K^R)^*f +\text{terms of order zero in}\, u.
\end{aligned}
\end{equation}

Next, we shall estimate
\begin{equation}
\label{1.26}
\partial_t\langle (K^R)^*u,(K^R)^*u\rangle=\frac{d\;}{dt} \int|(K^R)^*u|^2(x,t)dx.
\end{equation}

After  using the equation \eqref{1.25} to estimate  \eqref{1.26} we separate the terms obtained
into four groups $I_{K;j},\;j=1,..,4$.

In $I_{K;1}$ we set the terms with coefficients independent of $t$, i.e. those involving the operators
$$
\mathcal{L}^R(x),\,\,\,\vec b_1(x,0)\cdot\nabla,\,\,\,\vec b_2(x,0)\cdot\nabla,\,\,\,\mathcal{E}^R(x),
$$
in the equation.

In $I_{K;2}$ one has the terms involving the difference of the coefficients at the time $t$ and time $0$, i.e.
the terms containing the operators
$$
(\mathcal{L}(x,t)-\mathcal{L}(x,0)),\,(\vec b_1(x,t)-\vec b_1(x,0))\cdot\nabla,\,(\vec 
b_2(x,t)-\vec b_2(x,0))\cdot\nabla.
$$

$I_{K;3}$ contains the terms coming from the $\epsilon$ viscosity part of the equation, i.e.
$$
-\epsilon\langle (K^R)^*\Delta^2 u, (K^R)^*u\rangle -\epsilon\langle(K^R)^*, (K^R)^*\Delta^2u\rangle.
$$

In $I_{K;4}$ we collect the terms of order zero and those coming from the inhomogeneous term $f(x,t)$.
From the $L^2$-continuity of $(K^R)^*$ they are bounded by
\begin{equation}
\label{1.27}
cR^{2N_0}\|u(t)\|_2^2+cR^{2N_0}\|u(t)\|_2\|f\|_2,
\end{equation}
with  $N_0$ depending only on the dimension $n$.

Similarly, we shall estimate

\begin{equation}
\label{1.28}
\partial_t\langle E^Ru, E^Ru\rangle = \frac{d\;}{dt} \int|E^Ru|^2(x,t)dx.
\end{equation}
where
$$
I=E^R+\tilde K^R(K^R)^*,
$$
see Lemma 5.2.6 in \cite{KePoRoVe}.

Thus,
\begin{equation}
\label{1.29}
\begin{aligned}
&\partial_t E^Ru=-\epsilon E^R \Delta^2u +  iE^R \mathcal{L}(x,t)u+E^R\vec b_1(x,t)\cdot\nabla u\\
&\;\;\;\;+E^R\vec b_2(x,t)\cdot\nabla 
\bar u+E^Rc_1(x,t)u+E^Rc_2(x,t)\bar u +E^R f.
\end{aligned}
\end{equation}

 Inserting  the identity

\begin{equation}
\label{1.29b}
\begin{aligned}
&E^R \mathcal{L}(x,t)=E^R \mathcal{L}(x,0) + E^R (\mathcal{L}(x,t)-\mathcal{L}(x,0))\\
&=\mathcal{L}(x,0) E^R + [E^R;\mathcal{L}(x,0)]+(\mathcal{L}(x,t)-\mathcal{L}(x,0)) E^R\\
&+[ E^R;\mathcal{L}(x,t)-\mathcal{L}(x,0)],
\end{aligned}
\end{equation}

into the equation \eqref{1.29} and \eqref{1.28}
we split the terms obtained  into four groups 
$I_{E;j},\,j=1,..,4$.

In $I_{E;1}$ we place the terms with coefficients independent of $t$ coming from the expressions
$$
\mathcal{L}^R(x),\,\,\,\,\vec b_1(x,0)\cdot\nabla,\,\,\,\,\vec b_2(x,0)\cdot\nabla,\,\,\,\,\mathcal{E}^R(x),
$$
in the equation.

 $I_{E;2}$ contains the terms involving the operator
$$
(\mathcal{L}(x,t)-\mathcal{L}(x,0)),
$$
(in this case the terms involving the operators $E^R\Psi_{b_1},\;E^R\vec b_2(x,t)\cdot\nabla$ can be directly
bounded by using Lemma 5.2.6 in \cite{KePoRoVe}).

$I_{E;3}$ contains the terms coming from the $\epsilon$ viscosity part of the equation, i.e.
$$
-\epsilon\langle E^R\Delta^2 u, E^Ru\rangle -\epsilon\langle E^R, E^R\Delta^2u\rangle.
$$

In $I_{E;4}$ we collect the remainder terms of order zero 
which using Lemma 5.2.6 in \cite{KePoRoVe} are bounded by

\begin{equation}
\label{1.30g}
cR^{2N_0}\|u(t)\|_2^2+cR^{2N_0}\|u(t)\|_2\|f\|_2,
\end{equation}
with $N_0$ as in \eqref{1.27}.

 The terms in $I_{K;1}$ and $I_{E;1}$ (all involving time independent coefficients)
were already considered in  Section 5 of \cite{KePoRoVe}, (i.e. their contribution to the equations
of \eqref{1.26} and \eqref{1.28}). Thus, to bound them we shall use the following inequalities proved 
in Lemmas 4.2.4-5.2.6 in \cite{KePoRoVe}.

\begin{proposition}
\label{p0}
 There exist $M=M(N)\in\zz^+$, $N$ denoting the decay of the coefficients 
in \eqref{1.6}-\eqref{1.9.b}, with $M(N)\uparrow \infty$ as $N\uparrow \infty$, and 
$R_0$ large enough such that for 
$R\geq R_0$ 
one has 
\begin{equation} 
\label{1.30a}
\mathcal{R}e\,i\,\langle \mathcal{L}^R(x)(K^R)^*u,(K^R)^*u\rangle\equiv 0,
\end{equation}
\begin{equation}
\label{1.30b}
\mathcal{R}e\,i\,\langle \mathcal{L}(x,0)E^Ru,E^Ru\rangle\equiv 0,
\end{equation}
\begin{equation}
\label{1.30c}
\begin{aligned}
&|\langle i[\mathcal{L}^R;(K^R)^*]u+(K^R)^*\vec b_1(x,0)\cdot \nabla u,(K^R)^*u\rangle|\\
&+
|\langle \mathcal{E^R}(K^R)^*u,(K^R)^*u\rangle|+
|\langle [\mathcal{E}^R;(K^R)^*]u,(K^R)^*u\rangle|\\
&+|\langle (K^R)^*\vec b_2(x,0)\cdot \nabla \bar u ,(K^R)^*u\rangle|\\
&\leq c_0 R^{-M}\|(J^{1/2}u)\langle x\rangle^{-\tilde N}\|_2^2+c_0R^{N_0+N_1(M)}\|u(t)\|_2^2
+c_0\|f\|_2^2.
\end{aligned}
\end{equation}
and
\begin{equation}
\label{1.30d}
\begin{aligned}
&|\langle [E^R;\mathcal{L}(x,0)]u, E^Ru\rangle|\\
&+|\langle E^R\vec b_1(x,0)\cdot \nabla u , E^Ru\rangle|+
|\langle E^R\vec b_2(x,0)\cdot \nabla \bar u , E^Ru\rangle|\\
&\leq c_0R^{N_0}\|u(t)\|_2^2
+c_0\|f\|_2^2,
\end{aligned}
\end{equation}
with $N_0$ as in \eqref{1.27}, $\tilde N$ depending only on $n$ , with $N\gg \tilde N$, and $N_1$ depending only on 
$M$.
\end{proposition}

To handle the contribution coming from the terms in $I_{K;2}$
we shall establish the following inequalities.

\begin{proposition}
\label{p1}
 
 Take $N_0$ as in \eqref{1.27}, then
\begin{equation}
\label{1.32}
\begin{aligned}
&|\langle [(\mathcal{L}(x,t)-\mathcal{L}(x,0));(K^R)^*]u,(K^R)^*u\rangle|
\\
&\leq c T R^{N_0}\|(J^{1/2}u)\langle x\rangle^{-\tilde N}\|_2^2+c R^{N_0}\|u(t)\|_2^2,
\end{aligned}
\end{equation}
\begin{equation}
\label{1.33}
\mathcal{R}e\,i\,\langle (\mathcal{L}(x,t)-\mathcal{L}(x,0))(K^R)^*u,(K^R)^*u\rangle\equiv 0,
\end{equation}
\begin{equation}
\label{1.34}
\begin{aligned}
&|\langle (K^R)^*(\vec b_1(x,t)-\vec b_1(x,0))\cdot \nabla u ,(K^R)^*u\rangle|\\
&\leq c T R^{N_0}\|(J^{1/2}u)\langle x\rangle^{-\tilde N}\|_2^2+c R^{N_0}\|u(t)\|_2^2,
\end{aligned}
\end{equation}
and
\begin{equation}
\label{1.35}
\begin{aligned}
&|\langle (K^R)^*(\vec b_2(x,t)-\vec b_2(x,0))\cdot \nabla \bar u ,(K^R)^*\bar u\rangle|\\
&\leq c T R^{N_0}\|(J^{1/2}u)\langle x\rangle^{-\tilde N}\|_2^2+c R^{N_0}\|u(t)\|_2^2,
\end{aligned}
\end{equation}
with $\tilde N$ depending only on $n$ and $N\gg \tilde N$.
\end{proposition}

\begin{proof}
We write
$$
\aligned
&([(\mathcal{L}(x,t)-\mathcal{L}(x,0));(K^R)^*])^*\\
&=(\mathcal{L}(x,t)-\mathcal{L}(x,0))^*K^R-K^R(\mathcal{L}(x,t)-\mathcal{L}(x,0))^* 
\endaligned
$$
where
$$
\aligned
&(\mathcal{L}(x,t)-\mathcal{L}(x,0))^*\\
&=(a_{jk}(x,t)-a_{jk}(x,0))\partial^2_{jk}+\partial_j(a_{jk}(x,t)-a_{jk}(x,0))\partial_k\\
&=\beta_{jk}(x,t)\partial^2_{jk}+\tilde \beta_k(x,t)\partial_k=\Theta_1+\Theta_2,
\endaligned
$$
with the notation
$\partial_j$ instead of $\partial_{x_j}$.
Notice that the coefficients 
$\beta_{jk}(x,t)$ and $\tilde\beta_k(x,t)$ and their derivatives
up to order $p$ (large enough) have strong decay at infinity uniformly in $t\in[0,T]$.
 In fact, a sufficient number of  semi-norms in 
the Schwartz 
class are 
bounded by $cT$ for all $t\in[0,T]$.

The term $\Theta_2$ yields 
$$
\langle (\tilde \beta_k(x,t)\partial_k K^R-K^R\tilde \beta_k(x,t)\partial_k)u,(K^R)^*u\rangle
$$
which using Theorem 3.3.1 in \cite{KePoRoVe} and  the $L^2$-continuity of
$K^R$ can be bounded by
$$
T R^{N_0} \|u(t)\|_2^2,\;\;\;\;\;\;\;\forall t\in(0,T).
$$

 To handle $\Theta_1$ we use again Theorem 3.3.1 in \cite{KePoRoVe} to write that
\begin{equation}
\label{int. factor}
\beta_{jk}(x,t)\partial^2_{jk}K^R-K^R\beta_{jk}(x,t)\partial^2_{jk}=\Psi_{d^R}+\text{ zero order terms},
\end{equation}
with
$$
d^R(x,t,\xi)=
\beta_{jk}(x,t)\partial_{\xi_l}(\xi_j\xi_k)\partial_{x_l}K^R-
\partial_{x_l}\beta_{jk}(x,t)\xi_j\xi_k\partial_{\xi_l}K^R.
$$
We recall that a symbol in our class $a(x,\xi)$  when multiplied by $\phi(x)$, a fast decaying function
in the $x$ variable, becomes a classical symbol of the same order.

 To each term in the symbol of  $d^R(x,t,\xi)$ we can apply the following argument.

\underline{Claim} : Let $\phi(x,t)$ be an smooth function with strong decay at infinity in the 
$x$-variable uniformly in $t\in[0,T]$, with a sufficient number of  semi-norms in the Schwartz class bounded by $cT$. Then
\begin{equation}
\label{1.36}
|\langle \phi(x,t)\partial_{j}K^Ru,(K^R)^*u\rangle| \leq  
cTR^{2N_0}\|(J^{1/2}u(t))\langle x\rangle^{-\tilde N}\|_2^2+\|u(t)\|_2^2
\end{equation}

\underline {Notation} We shall use the notation $\simeq$ to denote terms whose difference can be bounded by a multiple
of $\|u\|_2^2$ or by an operator of order zero.

\begin{proof}
We have
$$
\aligned
&\langle \phi(x,t)\partial_{j}K^Ru,(K^R)^*u\rangle\\
&=\langle \langle x\rangle^{-\tilde N} \langle x\rangle^{2\tilde N} \phi(x,t)\partial_{j}K^RJ^{1/2}J^{-1}J^{1/2}u, \langle 
x\rangle^{-\tilde N}(K^R)^*u\rangle\\
&\simeq \langle \langle x\rangle^{-\tilde N} J^{1/2}\langle x\rangle^{2\tilde N} \phi(x,t)\partial_{j}K^R J^{-1}J^{1/2}u, \langle 
x\rangle^{-\tilde N}(K^R)^*u\rangle\\
&\simeq \langle \langle x\rangle^{-\tilde N} J^{1/2}\langle x\rangle^{3\tilde N} \phi(x,t)\partial_{j}K^R J^{-1}\langle 
x\rangle^{-\tilde N}J^{1/2}u, 
\langle
x\rangle^{-\tilde N}(K^R)^*u\rangle\\
&\simeq \langle J^{1/2}\langle x\rangle^{-\tilde N} \langle x\rangle^{3 \tilde N} \phi(x,t)\partial_{j}K^R J^{-1}\langle 
x\rangle^{-\tilde N}J^{1/2}u,
\langle
x\rangle^{-\tilde N}(K^R)^*u\rangle\\
&\simeq\langle\langle x\rangle^{3\tilde N} \phi(x,t)\partial_{j}K^R J^{-1}\langle x\rangle^{-\tilde N}J^{1/2}u,
\langle x\rangle^{-\tilde N}J^{1/2}\langle
x\rangle^{-\tilde N}(K^R)^*u\rangle\\
&=\Lambda_3.
\endaligned
$$
 Since $\langle x\rangle^{3\tilde N} \phi(x,t)\partial_{j}K^R J^{-1}\in S^0_{1,0}$, i.e. a classical $\psi$.d.o.
of order zero, it follows that
$$
|\Lambda_3|\leq cTR^{N_0}\|\langle x\rangle^{-\tilde N}J^{1/2}u\|_2\|\langle x\rangle^{-\tilde N}J^{1/2}\langle 
x\rangle^{-\tilde N}(K^R)^*u\|_2.
$$
Finally, using that
$$
(\langle x\rangle^{-\tilde N}J^{1/2}\langle x\rangle^{-\tilde N}(K^R)^*)^*=
K^R\langle x\rangle^{-\tilde N}J^{1/2}\langle x\rangle^{-\tilde N},
$$
we obtain the desired inequality \eqref{1.36}.

 Returning to the operator $\Psi_{d^R}$ whose symbol is $d^R(x,\xi)$ and applying our claim one gets that
$$
|\langle \Psi_{d^R} u, (K^R)^*u\rangle|\leq cTR^{N_0}\|\langle x\rangle^{-\tilde N}J^{1/2}u(t)\|_2^2+\|u(t)\|_2^2,
$$
which proves \eqref{1.32}.

 Integration by parts yields  \eqref{1.33}.

To prove  \eqref{1.34}  we write
$$
((K^R)^*(b_1(x,t)-b_1(x,0))\cdot\nabla))^*\simeq ((b_1(x,t)-b_1(x,0))\cdot\nabla)K^R=G_R^*\in S^1_{1,0},
$$
and $\Psi_{G_R^*}\simeq \Psi_{\bar G_R}$, so using an argument
similar to that given in the proof of the claim above we get  \eqref{1.34}. Similarly for  \eqref{1.35}.

\end{proof}

 To bound the terms in $I_{E;2}$ we have the following estimates.

\begin{proposition}
\label{p2}
There exists $N_0\in\zz^+$ depending only on the dimension $n$ such that for any $t\in (0,T_0]$
\begin{equation}
\label{1.37}
|\langle [(\mathcal{L}(x,t)-\mathcal{L}(x,0));E^R]u,E^Ru\rangle|
\leq c TR^{N_0}\|u(t)\|_2^2,
\end{equation}
and
\begin{equation}
\label{1.37b}
\mathcal{R}e \,i\langle (\mathcal{L}(x,t)-\mathcal{L}(x,0))E^Ru,E^Ru\rangle=0,
\end{equation}

\end{proposition}

\begin{proof}
 
With the same notation that in the previous proof we have
$$
(\mathcal{L}(x,t)-\mathcal{L}(x,0))u=\partial_j(\beta_{jk}(x,t)\partial_k u).
$$
So using that $E^R=I-\tilde K^R(K^R)^*$, one sees that
$$
\aligned
& [(\mathcal{L}(x,t)-\mathcal{L}(x,0));E^R]u \\
&=\tilde K^R[(K^R)^*;\partial_j\beta_{jk}(x,t)\partial_k]+[\tilde K^R;\partial_j\beta_{jk}(x,t)\partial_k](K^R)^*.
\endaligned
$$
First we consider the term involving $[(K^R)^*;\partial_j\beta_{jk}(x,t)\partial_k]$. 
Taking the adjoint it follows that
$$
\aligned
&([(K^R)^*;\partial_j\beta_{jk}(x,t)\partial_k])^*=
- K^R\partial_j\beta_{jk}(x,t)\partial_k+\partial_j\beta_{jk}(x,t)\partial_k  K^R\\
&=- K^R\partial_j(\beta_{jk}(x,t))\partial_k+\partial_j(\beta_{jk}(x,t))\partial_k  K^R\\
&\;- K^R\beta_{jk}(x,t)\partial^2_{jk}+\beta_{jk}(x,t)\partial^2_{jk} K^R\\
&\simeq - K^R\beta_{jk}(x,t)\partial^2_{jk}+\beta_{jk}(x,t)\partial^2_{jk} K^R,
\endaligned
$$
since $- K^R\partial_j(\beta_{jk})\partial_k+\partial_j(\beta_{jk})\partial_k  K^R$ is an operator
of order zero, see Theorem 3.3.1 in \cite{KePoRoVe}.

 Up to symbols of order zero (bounded operator in $L^2$) the symbol of $\beta_{jk}(x,t)\partial^2_{jk} K^R- 
K^R\beta_{jk}(x,t)\partial^2_{jk}$  is equal to
$$
\eta(x,t,\xi)=\beta_{jk}\partial_{\xi_l}(\xi_k\xi_j)\partial_{x_l}k^R(x,\xi)-\partial_{x_l}\beta_{jk}\xi_k\xi_j
\partial_{\xi_l}k^R(x,\xi).
$$
 So $\eta(x,t,\xi)\in S^1_{1,0}$ uniformly in $t\in[0,T]$ and $ \eta(x,t,\xi)\simeq \Psi_{d}\tilde 
\beta_l(x,t)\partial_l$
with $d=d(x,\xi)$ in our class. Similarly for $\eta^*(x,t,\xi)$, then
$\eta^*(x,t,\xi)\simeq \tilde \beta_l(x,t)\partial_l \Psi_{d_1}$. Inserting this in \eqref{1.37}
we get
$$
\aligned
&\langle \tilde K^R [(K^R)^*;\partial_j\beta_{jk}(x,t)\partial_k]u,E^Ru\rangle
\simeq \langle \tilde K^R \tilde \beta_l(x,t)\partial_l \Psi_{d_1(x,\xi)}u,E^Ru\rangle\\
&\simeq \langle \Psi_b \Psi_{d_1(x,\xi)}u,E^Ru\rangle \simeq
\langle \Psi_{d_1(x,\xi)}u, ( \Psi_b)^*E^Ru\rangle,
\endaligned
$$
which can be bounded by $\|u(t)\|^2_2$, by using that $\Psi_b$ is a classical $\psi$.d.o. of order 1
and the continuities properties of $E^R$ deduce in Lemma 5.2.6 in \cite{KePoRoVe}.
 
A similar argument provides the bound for the term $[\tilde K^R,\partial_j\beta_{jk}(x,t)\partial_k](K^R)^*$.
Collecting this information we get \eqref{1.37}.

The proof of \eqref{1.37b} follows by integration by parts.

\end{proof}

 The terms coming from the artificial viscosity term $\epsilon \Delta^2$, i.e.
the terms in $I_{K;3}$ and $I_{E;3}$ will be handled by using the following 
inequalities. 

\begin{proposition}
\label{p3}

\begin{equation}
\label{1.37c}
\langle \epsilon (K^R)^*\Delta^2u, (K^R)^*u\rangle = \epsilon \langle \Delta(K^R)^*u,\Delta(K^R)^*u\rangle + \Lambda_1,
\end{equation}
and
\begin{equation}
\label{1.37f}
\langle \epsilon E^R\Delta^2u,  E^Ru\rangle = \epsilon \langle \Delta  E^R u,\Delta  E^R u\rangle + \Lambda_2,
\end{equation}
where
$$
|\Lambda_j|\leq \epsilon R^{2N_0}\|\Delta u(t)\|_2(\|u(t)\|_2+\|\nabla u(t)\|_2),\;\;\;j=1,2.
$$
\end{proposition}

\begin{proof}

To obtain \eqref{1.37c} we write
\begin{equation}
\label{1.37d}
\begin{aligned}
&\langle (K^R)^*\Delta^2u, (K^R)^*u\rangle\\
&=\langle [(K^R)^*;\Delta] \Delta u, (K^R)^*u\rangle
+\langle \Delta (K^R)^*\Delta u, (K^R)^*u\rangle\\
&=\langle [(K^R)^*;\Delta] \Delta u, (K^R)^*u\rangle
+\langle (K^R)^*\Delta u, \Delta (K^R)^*u\rangle\\
&=\langle [(K^R)^*;\Delta] \Delta u, (K^R)^*u\rangle
+\langle(K^R)^* \Delta u,  [\Delta; (K^R)^*]u\rangle \\
&+\langle(K^R)^* \Delta u,(K^R)^* \Delta u\rangle\\
&\equiv  \Omega_1+\Omega_2+\Omega_3.
\end{aligned}
\end{equation}

 Thus, we consider
$$
 [\Delta; (K^R)^*]=\Delta (K^R)^*-(K^R)^*\Delta =\mathcal{T},
$$
and its adjoint
$$
\mathcal{T}^*=-(\Delta K^R - K^R\Delta).
$$
Let $k^R(x,\xi)$ denote the symbol of $ K^R$
such that
$$
K^Rf(x)=\int e^{ix\cdot \xi}k^R(x,\xi)\hat f(\xi)d\xi,
$$
so
$$
\Delta K^Rf(x)=\int e^{ix\cdot \xi}\{-|\xi|^2k^R(x,\xi)+2i\xi_j\partial_{x_j}k^R(x,\xi)+\Delta_x k^R(x,\xi)\}\hat 
f(\xi)d\xi.
$$
Therefore,
$$
\mathcal{T}^*f(x)=\int e^{ix\cdot \xi}\{2i\xi_j\partial_{x_j}k^R(x,\xi)+\Delta_xk^R(x,\xi)\}\hat f(\xi)d\xi.
$$

 From Lemma \ref{ccd} above  and Theorem 3.2.1 in \cite{KePoRoVe} one has that the operators with symbols
$\partial_{x_j}k^R(x,\xi)$ and $\Delta_xk^R(x,\xi)$ are bounded in $L^2$. Hence, we can write that
\begin{equation}
\label{number1}
\mathcal{T}^*=C_j\partial_{x_j}+C_0,\;\;\;\text{and}\;\;\;\mathcal{T}=C_0^*+\partial_{x_j}(-C_j^*),
\end{equation}
with $C_j,\;j=0,1,..,n$ denoting $L^2$-bounded operators.

Also we have that
$$
\aligned
&\int e^{ix\cdot \xi}i\xi_j\partial_{x_j}k^R(x,\xi)\hat f(\xi)d\xi\\
&=\partial_{x_j}(\int e^{ix\cdot \xi}\partial_{x_j}k^R(x,\xi)\hat f(\xi)d\xi)-\int e^{ix\cdot\xi}
\partial^2_{x_jx_j}k^R(x,\xi)\hat f(\xi)d\xi,
\endaligned
$$
so that
\begin{equation}
\label{number2}
\mathcal{T}^*=\partial_{x_j}\tilde C_j+\tilde C_0,
\;\;\;\text{and}\;\;\;\mathcal{T}=\tilde C_0^*-\tilde C_j^*\partial_{x_j},
\end{equation}
with $\tilde C_j,\;j=0,1,..,n$  denoting $L^2$ bounded operators.

To bound $\Omega_1$ we see that
$$
\aligned
&\langle  [(K^R)^*;\Delta]\Delta u,(K^R)^*u\rangle =-\langle \mathcal{T}\Delta u,(K^R)^*u\rangle\\
&=-\langle(C_0^*+\partial_{x_j}(-C_j^*))\Delta u,(K^R)^*u\rangle \\
&= 
-\langle C_0^*\Delta u,(K^R)^*u\rangle - \langle (-C_j^*)\Delta u,\partial_{x_j}(K^R)^*u\rangle.
\endaligned
$$
Since an explicit computation shows that
$$
\partial_{x_j}(K^R)^*\simeq (K^R)^*\partial_{x_j},
$$
i.e. their difference is an $L^2$-bounded operator, it follows that
$$
\aligned
&|\Omega_1|\leq c_0\|\Delta u\|_2(\|u\|_2+\|(K^R)^*\nabla u\|_2)\\
&\leq c_0\|\Delta u\|_2(\|u\|_2+R^{N_0}\|\nabla u\|_2).
\endaligned
$$

Using again \eqref{number1}-\eqref{number2} we have that
$$
\aligned
&\langle (K^R)^*\Delta u, [\Delta;(K^R)^*]u\rangle=\langle (K^R)^*\Delta u,\mathcal{T} u\rangle
&=\langle (K^R)^*\Delta u,(\tilde C_0^*-\tilde C_j^*\partial_{x_j})u\rangle,
\endaligned
$$
so $\Omega_2$ can be bounded as
$$
\aligned
&|\Omega_2|\leq c_0\|(K^R)^* \Delta u\|_2(\|u\|_2 +\|\nabla u\|_2)\\
&\leq  c_0 R^{N_0} \|\Delta u\|_2(\|u\|_2+\|\nabla u\|_2).
\endaligned
$$

Inserting this information in \eqref{1.37d} we obtain \eqref{1.37c}.

 To obtain \eqref{1.37f} we recall that $E^R=I-\tilde K^R(K^R)^*$ so that
$$
E^R \Delta^2=\Delta E^R \Delta + [E^R;\Delta]\Delta,
$$
with
$$
\aligned
&[E^R;\Delta]=-[\tilde K^R(K^R)^*;\Delta]=-(\tilde K^R(K^R)^*\Delta -\Delta \tilde K^R(K^R)^*)\\
&=-\tilde K^R [(K^R)^*;\Delta]+[\tilde K^R;\Delta](K^R)^*\equiv \Gamma_1+\Gamma_2.
\endaligned
$$
 Using \eqref{number1}-\eqref{number2} and an explicit computation it follows that
$$
\Gamma_1=\tilde K^R \mathcal{T}=\tilde K^R( C_0^*-\partial_{x_j}C_j^*)= -\partial_{x_j}\tilde K^RC_j^*+Q_1,
$$
where $Q_j,\;j=1,..,4$ will denote $L^2$ bounded operators. So
$$
\aligned
&|\langle \Gamma_1\Delta u, E^Ru\rangle|= |\langle (\partial_{x_j}\tilde K^R C_j^*+Q_1)\Delta u, E^Ru\rangle|\\
&\leq c_0 \|\Delta u\|_2(\|\partial_{x_j} E^Ru\|_2+\| E^Ru\|_2).
\endaligned
$$
Combining  $\partial_{x_j}E^R=\partial_{x_j}-\partial_{x_j}\tilde K^R(K^R)^*$ and an explicit computation one gets 
that
$$
\aligned
&\|\partial_{x_j}E^Ru\|_2 = \|(\partial_{x_j}-\partial_{x_j}\tilde K^R(K^R)^*)u\|_2\\
&\leq \|u\|_2 + \|\partial_{x_j}\tilde K^R(K^R)^* u\|_2
\simeq \|u\|_2+\|\tilde K^R\partial_{x_j}(K^R)^* u\|_2\\ 
&\leq \|u\|_2 + R^{N_0}\|\partial_{x_j}(K^R)^* u\|_2
=\|u\|_2+ R^{N_0} \|K^R \partial_{x_j}  u\|_2\\
&\leq \|u\|_2 +R^{2N_0} \| \partial_{x_j}  u\|_2,
\endaligned
$$
and consequently
$$
|\langle\Gamma_1\Delta u,E^Ru\rangle|\leq c_0 R^{2N_0} \|\Delta u\|_2(\|u\|_2+\|\nabla u\|_2).
$$

 For $\Gamma_2$ we reproduce the argument in \eqref{number1}-\eqref{number2} for $\tilde K^R$ 
instead of $ K^R$. So using the same notation we have
$$
\Gamma_2=[\tilde K^R;\Delta](K^R)^*=(C_0+\partial_{x_j}C_j)(K^R)^*,
$$
so as before
$$
\aligned
&|\langle\Gamma_2\Delta u,E^Ru\rangle|=
|\langle (C_0+\partial_{x_j}C_j)(K^R)^*\Delta u,E^Ru\rangle|\\
&\leq |\langle C_0(K^R)^*\Delta u,E^Ru\rangle|+|\langle C_j(K^R)^*\Delta u,\partial_{x_j}E^Ru\rangle|\\
&\leq \|\Delta u\|_2(\|E^Ru\|_2+\|\partial_{x_j}E^Ru\|_2)
\leq c_0 R^{2N_0}\|\Delta u\|_2(\|u\|_2+\|\nabla u\|_2).
\endaligned
$$

 Collecting the results in Propositions \ref{p0}-\ref{p3}
we get that
\begin{equation}
\label{1.38}
\begin{aligned}
&\frac{d\;}{dt}\|(K^R)^*u(t)\|_2^2 +\epsilon\|(K^R)^*\Delta u(t)\|_2^2\\
&\leq (c_0 R^{-M}+cTR^{N_0})\|(J^{1/2}u)(t)\langle x\rangle^{-\tilde N}\|_2^2 +  c 
(R^{N_0+N_1(M)}+TR^{N_0})\|u(t)\|_2^2\\
&+c_0\|f\|_2^2 
+c_0\epsilon R^{2N_0}\|\Delta u(t)\|_2(\|u(t)\|_2+\|\nabla u(t)\|_2),
\end{aligned}
\end{equation}
and
\begin{equation}
\label{1.39}
\begin{aligned}
&\frac{d\;}{dt}\|E^Ru\|_2^2 +\epsilon\|E^R\Delta u\|_2^2\\
& \leq c (R^{N_0}+TR^{N_0})\|u(t)\|_2^2
+c\|f\|_2^2 \\
&+c_0\epsilon R^{2N_0}\|\Delta u(t)\|_2(\|u(t)\|_2+\|\nabla u(t)\|_2).
\end{aligned}
\end{equation}

We will use that $E^R=I-\tilde K^R(K^R)^*$ so that
\begin{equation}
\label{1.40}
\begin{aligned}
&\|v\|_2=\|(E^R+\tilde K^R(K^R)^*)v\|_2\\
&\leq \|E^R v\|_2 +  R^{N_0}\|(K^R)^*v\|_2.
\end{aligned}
\end{equation}
Also, we need the following interpolation estimates : for any $\,l\in\zz^+$
\begin{equation}
\label{1.41}
R^l\|v\|_2 \|\Delta v(t)\|_2\leq R^{2l}\|v\|_2^2 +\|\Delta v\|_2^2,
\end{equation}
and
\begin{equation}
\label{1.42}
R^l\|\nabla v\|_2 \|\Delta v\|_2
\leq cR^l\|v\|_2^{1/2} \|\Delta v\|_2^{3/2} 
\leq c R^{4l}\|v\|_2^2 +\|\Delta v\|_2^2.
\end{equation}

So combining \eqref{1.38}-\eqref{1.39} with \eqref{1.40}-\eqref{1.42}
we find that
\begin{equation}
\label{1.43}
\begin{aligned}
&\frac{d\;}{dt}(\|(K^R)^*u(t)\|_2^2+\|E^Ru(t)\|^2_2) +\epsilon R^{-2N_0}\|\Delta u(t)\|_2^2\\
&\leq (cTR^{N_0}+c_0R^{-M})\|(J^{1/2} u)\langle x\rangle^{-\tilde N}\|_2^2 + c(R^{N_0+N_1(M)}+TR^{N_0})
\|u(t)\|_2^2 \\
&+ c_0\|f(t)\|_2^2 + c_0\epsilon R^{2N_0}\|\Delta u(t)\|_2(\|u(t)\|_2+\|\nabla u(t)\|_2)\\
&\leq (cTR^{N_0}+c_0R^{-M})\|(J^{1/2} u)\langle x\rangle^{-\tilde N}\|_2^2 + c(R^{N_0+N_1(M)}+TR^{N_0})
\|u(t)\|_2^2 \\
&+ c_0\|f(t)\|_2^2 +\frac{\epsilon}{2}R^{-2N_0}\|\Delta u(t)\|_2^2 +c_0\epsilon R^{14N_0}\|u(t)\|_2^2.
\end{aligned}
\end{equation}
Thus,
\begin{equation}
\label{1.44}
\begin{aligned}
&\frac{d\;}{dt}(\|(K^R)^*u(t)\|_2^2+\|E^Ru(t)\|^2_2) +\frac{\epsilon}{2} R^{-2N_0}\|\Delta u(t)\|_2^2\\
&\leq (cTR^{N_0}+c_0R^{-M})\|(J^{1/2} u)\langle x\rangle^{-\tilde N}\|_2^2 + c(R^{N_0+N_1(M)}+TR^{N_0})
\|u(t)\|_2^2 \\
&+ c_0\|f(t)\|_2^2+c_0\epsilon R^{14N_0}\|u(t)\|_2^2.
\end{aligned}
\end{equation}

Integrating the above inequality in the time interval $(0,T)$ and inserting in the result
the estimate \eqref{1.19} we get
\begin{equation}
\label{1.45}
\begin{aligned}
&\sup_{0\leq t\leq T} (\|(K^R)^*u(t)\|_2^2+\|E^Ru(t)\|^2_2)
+\frac{\epsilon}{2} R^{-2N_0}\int_0^T\|\Delta u(t)\|_2^2dt\\
&\leq \|(K^R)^*u(0)\|_2^2+\|E^Ru(0)\|^2_2\\
&+(cTR^{N_0}+c_0R^{-M})\int_0^T\|(J^{1/2} u)\langle x\rangle^{-\tilde N}\|_2^2dt\\
&+(c R^{N_0+N_1(M)}+TR^{N_0}+c_0\epsilon R^{14N_0})\int_0^T \|u(t)\|_2^2dt+c_0\int_0^T\|f(t)\|_2^2dt\\
&\leq \|(K^R)^*u(0)\|_2^2+\|E^Ru(0)\|^2_2\\
&+(cTR^{N_0}+c_0R^{-M})((c_0+cT)\sup_{0\leq t\leq T}\|u(t)\|_2^2\\
&+c_0\epsilon\int^T_0 \|\Delta u(t)\|^2_2dt +c_0\int_0^T\|f(t)\|_2^2dt)\\
&+c (R^{N_0+N_1(M)}+TR^{N_0}+c_0\epsilon R^{14N_0})T\sup_{0\leq t\leq T} \|u(t)\|_2^2+c_0\int_0^T\|f(t)\|_2^2dt.
\end{aligned}
\end{equation}

Since $M(N)\uparrow \infty$ as $N\uparrow \infty$, we take $N$ in our hypotheses large enough 
such that $M=100N_0$ (we recall that $N_0$ depends only on the dimension $n$). Next,  
we fix $R\geq R_0$ sufficiently large and then $T=T(N_0,M,R)>0$ small enough 
such that the following inequalities holds
$$
(cTR^{N_0}+c_0R^{-M})(c_0+cT)\leq \frac{R^{-2N_0}}{4},
$$
$$
(cR^{N_0+N_1(M)}+cTR^{N_0}+c_0 
R^{14N_0})T\leq  \frac{R^{-2N_0}}{4}
$$
and
$$
c_0(cTR^{N_0}+c_0R^{-M})\leq R^{-N_0}/4.
$$

 Combining these inequalities, \eqref{1.45}, and \eqref{1.40} we get the estimate
$$
\sup_{0\leq t\leq T}\|u(t)\|_2^2+
\frac{\epsilon}{4}R^{-2N_0}\int_0^T\|\Delta u(t)\|^2_2dt\leq c_0R^{2N_0}\|u(0)\|_2^2+c_0\int_0^T\|f(t)\|_2^2dt.
$$
which proves \eqref{1.12}.

\end{proof} 
\end{proof}

\vspace{3mm} \noindent{\large {\bf Acknowledgments}}
\vspace{3mm}\\  C. E. K. and G. P. were supported by NSF grants. L. V. was supported by a DGICYT grant.
This work was carried out at the IAS which support we acknowledge.
C. E. K. was supported at IAS by The Von Newmann Fund, The Weyl
Fund, The Oswald Veblem Fund and The Bell Companies Fellowship.


\vskip1cm
\noindent{\bf Carlos E. Kenig}\\
Department of Mathematics\\
University of Chicago\\
Chicago, Il. 60637 \\
USA\\
E-mail: cek@math.uchicago.edu\\

\noindent{\bf Gustavo Ponce}\\
Department of Mathematics\\
University of California\\
Santa Barbara, CA 93106\\
USA\\
E-mail: ponce@math.ucsb.edu\\

\noindent{\bf Christian  Rolvung}\\
Nykredit Markets \& Asset Management\\
Kalvevod Bridge 1-3\\
DK-1780 Copenhagen V\\
Denmark\\
E-mail: cro@nykredit.dk\\

\noindent{\bf Luis Vega}\\
Departamento de Matematicas\\
Universidad del Pais Vasco\\
Apartado 644\\
48080 Bilbao\\
Spain\\
E-mail: mtpvegol@lg.ehu.es


\begin{thebibliography}{99}
\small

\bibitem{AbHa}  M. J. Ablowitz, R. Haberman,
{\em Nonlinear evolution equations in two and three dimensions},
Phys. Rev. Lett. {\bf 35} (1975), 1185--1188.


\bibitem{BoSm} J. L. Bona, R. Smith,
{\em  The initial value problem for the Korteweg-de Vries equation},
 Roy. Soc. London {\bf Ser A 278} (1978),  555--601.




\bibitem{CaVa} A. P. Calder\'on, R. Vaillancourt,
{\em On the boundedness of pseudo-differential operators},
J. Math. Soc. Japan
{\bf 23} (1971), 374--378.



\bibitem{c-s} P. Constantin, J-C. Saut,
{\em Local smoothing properties of dispersive equations}, Journal  A.M.S.
{\bf 1} (1989), 413--446.


\bibitem{Ch} H. Chihara,
{\em Local existence for semi-linear Schr\"odinger equations}, Math. Japan
{\bf 42} (1995), 35--51.

\bibitem{CrKaSt} W.  Craig, T. Kappeler, W. Strauss, 
{\em Microlocal dispersive smoothing for the Schr\"odinger equation},
 Comm. Pure Appl. Math.
{\bf  48} (1995), 769--860.

\bibitem{DaSt} A. Davey, K. Stewartson,
{\em On three--dimensional packets of surface waves},
 Proc. Roy. Soc. London Ser.
{\bf 338} (19740), 101--110.

\bibitem{DjRe}  V. D. Djordjevic,  L. G. Redekopp,
{\em  On two-dimensional packets of capillary-gravity waves}
J. Fluid Mech.
{\bf 79} (1977), 703--714.

\bibitem{Do1} S. Doi, 
{\em On the Cauchy problem for Schr\"odinger type equations and the regularity of solution}
 J. Math. Kyoto Univ.
{\bf 34} (1994), 319--328.

\bibitem{Do2} S. Doi, 
{\em Remarks on the Cauchy problem for Schr\"odinger--type equations}
 Comm. Partial Differential Equations
{\bf 21} (1996), 163--178.

\bibitem{Do3} S. Doi,
{\em  Smoothing effects for Schrodinger evolution equation and global behavior of geodesic flow}
 Math. Ann.
{\bf 318} (2000), 355--389.

\bibitem{HaOz}  N. Hayashi, T. Ozawa,
{\em Remarks on nonlinear Schr\"odinger equations in one space dimension}
Differential Integral Equations {\bf 7} (1994),  453--461.
 

\bibitem{ Ho}  L. H\"ormander,
{\em Pseudo--differential operators and non--elliptic boundary problems},
Ann. of Math.  {\bf 2} (1966),  129--209.

\bibitem{Ic} W. Ichinose, 
{\em On $L^2$ wellposedness of the Cauchy problem for Schr\"odinger 
type equations on a Riemannian manifold
and Maslov theory}, Duke Math. J
{\bf 56} (1988), 549--588.

\bibitem{Is} Y. Ishimori,
{\em Multi vortex solutions of a two dimensional nonlinear wave equation},
Progr. Theor. Phys {\bf 72} (1984), 33--37.

\bibitem{Kt} T. Kato,
 {\em On the Cauchy problem for the (generalized) Korteweg-de Vries equation},
Advances in Math. Supp. Studies, 
Studies in Applied Math.
{\bf 8} (1983), 93--128.

\bibitem{KePoVe1}  C. E. Kenig, G. Ponce, L. Vega, 
{\em Oscillatory integrals and regularity of dispersive equations},
Indiana University Math. J.
{\bf 40} (1991), 33--69.

\bibitem{KePoVe2}  C. E. Kenig, G. Ponce, L. Vega,
{\em Small solutions to nonlinear Schr\"odinger equations},
 Ann. Inst. H. Poincar\'e Anal. Non Lin\'eaire
{\bf 10} (1993), 255--288.

\bibitem{KePoVe3}  C. E. Kenig, G. Ponce, L. Vega,
{\em Smoothing effects and local existence theory for the generalized 
nonlinear Schr\"odinger equations}, 
Invent. Math. {\bf 134} (1998), 489-545.


\bibitem{KePoVe4}  C. E. Kenig, G. Ponce, L. Vega,
{\em The Cauchy problem for the quasi-linear Schr\"odinger equations},
Invent. Math. {\bf 158} (2004), 343-388.


\bibitem{KePoRoVe}  C. E. Kenig, G. Ponce, C. Rolvung, L. Vega, 
{\em Variable coefficient Schr\"odinger flows for ultrahyperbolic operators}
to appear in Advances in Math.


\bibitem{KrFa}  D. J. Kruzhkov, A. V. Faminskii,
{\em Generalized solutions for the Cauchy problem for the Korteweg-de Vries equation},
Math. USSR Sb, {\bf 48} (1990), 93--138.

\bibitem{Mi} S. Mizohata, 
{\em  On the Cauchy Problem}
Notes and Reports in Mathematics in Science and Engineering,3, 
Science Press and Academic Press
(1985).

\bibitem{Ro} C. Rolvung,
{\em Non-isotropic Schr\"odinger equations}
 PhD. dissertation, University of Chicago,
(1998).


\bibitem{s} P. Sj\"olin, {\em Regularity of solutions to the
Schr\"odinger equation}, Duke Math. J. {\bf 55} (1987), 699--715.

\bibitem{St}   E. M. Stein,
{\em Harmonic analysis: real--variable methods, orthogonality, and oscillatory integrals},
 Princeton University Press, (1993).

\bibitem{s-s} C. Sulem,  P-L. Sulem, {\em The nonlinear Schr\"odinger
Equation. Self-focusing and collapse}, Applied Mathematical Sciences {\bf 139}.
Springer (1999).

\bibitem{Ta} J. Takeuchi, {\em Le probl\`eme de Cauchy pour certaines equations aux deriv\'ees partielles du type de 
Schr\"odinger sym\'etrisations ind\'ependendantes du temps},  C. R. Acad. Sci. Paris, {\bf t315} serie 1 (1992),  1055-1058.




\bibitem{v} L. Vega, {\em Schr\"odinger equations: pointwise convergence
to the initial data},  Proc A.M.S. {\bf 102} (1988), 874--878.

\bibitem{ZaSc} V. E.  Zakharov,  E. I. Schulman, 
{\em Degenerated dispersion laws, motion invariant and kinetic equations},
 Physica {\bf 1D} (1980), 185-250.


\end{thebibliography}
\end{document}